\documentclass[10pt,twosided]{article}
\usepackage[tbtags]{amsmath}
\usepackage{amssymb,amsthm,latexsym}
\allowdisplaybreaks[4]
\usepackage{color}

\usepackage{amssymb,color}

\definecolor{c20}{rgb}{0.,0.7,0.}
\definecolor{c30}{rgb}{0.,0.,1.}
\definecolor{c40}{rgb}{1,0.1,0.7}
\definecolor{c50}{rgb}{1,0,0}
\definecolor{c60}{rgb}{1,0.9,0.1}

\newcommand{\EE}[1]{\mathbb{E}\left(#1\right)}
\def\qE#1{{\textcolor{c20}{#1}}}
\def\qE#1{#1}

\def\tE#1{{\textcolor{c20}{#1}}}
\def\tE#1{#1}

\def\tw#1{{\textcolor{c40}{#1}}}
\def\tw#1{#1}

\def\Ke#1{{\textcolor{c20}{#1}}}
\def\Ke#1{#1}

\def\Ee#1{{\textcolor{c40}{#1}}}
\def\Ee#1{#1}

\def\cW#1{{\textcolor{c40}{#1}}}
\def\cW#1{#1}
\def\wz#1{{\textcolor{c40}{#1}}}
\def\wz#1{#1}
\def\wzc#1{{\textcolor{c40}{#1}}}
\def\wzc#1{#1}

\def\He#1{\textcolor{c20}{#1}}
\def\He#1{#1}

\def\cE#1{\textcolor{c20}{#1}}
\def\cE#1{#1}

\def\cH#1{\textcolor{c20}{#1}}
\def\cH#1{#1}

\def\rH#1{\textcolor{c20}{#1}}
\def\rH#1{#1}

\def\aH#1{\textcolor{c20}{#1}}
\def\aH#1{#1}
\def\AH#1{\textcolor{c20}{#1}}
\def\AH#1{#1}

\def\zw#1{{\textcolor{c40}{#1}}}
\def\zw#1{#1}
\def\zc#1{{\textcolor{c40}{#1}}}
\def\zc#1{#1}

\def\zcw#1{{\textcolor{c40}{#1}}}
\def\zcw#1{#1}

\usepackage{amssymb,color}

\newcommand{\kb}[1]{\boldsymbol{#1}}
\newcommand{\vk}[1]{\kb{#1}}

\newcommand{\abs}[1]{\lvert #1 \rvert}

\newcommand{\E}[1]{\mathbb{E}\left(#1\right)}

\newcommand{\pk}[1]{\mathbb{P} \left( #1 \right) }

\newcommand{\EXP}[1]{\exp \left( #1 \right) }

\newcommand{\R}{\mathbb{R}}

\newcommand{\inr}{\in \R}

\newcommand{\ldot}{,\ldots,}

\newcommand{\limit}[1]{\lim_{#1 \to   \infty}}

\newcommand{\BQN}{\begin{eqnarray}}
\newcommand{\EQN}{\end{eqnarray}}
\newcommand{\BQNY}{\begin{eqnarray*}}
\newcommand{\EQNY}{\end{eqnarray*}}

\newcommand{\BS}{\begin{sat}}
\newcommand{\ES}{\end{sat}}
\newcommand{\BT}{\begin{theo}}
\newcommand{\ET}{\end{theo}}
\newcommand{\BK}{\begin{korr}}
\newcommand{\EK}{\end{korr}}

\newcommand{\BD}{\begin{de}}
\newcommand{\ED}{\end{de}}

\newcommand{\BIT}{\begin{itemize}}
\newcommand{\EIT}{\end{itemize}}
\newcommand{\BDI}{\begin{description}}
\newcommand{\EDI}{\end{description}}

\newcommand{\BRM}{\begin{remarks}}
\newcommand{\ERM}{\end{remarks}}

\newcommand{\BEL}{\begin{lem}}
\newcommand{\EEL}{\end{lem}}

\numberwithin{equation}{section}
\newtheorem{theo}{Theorem}[section]
\newtheorem{sat}[theo]{Proposition}
\newtheorem{de}[theo]{Definition}
\newtheorem{lem}{Lemma}[section]

\newtheorem{korr}[theo]{Corollary}
\newtheorem{remark}[theo]{Remark}
\newtheorem{remarks}[theo]{Remarks}

\newtheorem{theorem}{Theorem}[section]

\newtheorem{lemma}{Lemma}[section]

\newcommand{\prooftheo}[1]{ \textsc{Proof of Theorem} \ref{#1} }

\newcommand{\prooflem}[1]{\textsc{Proof of Lemma} \ref{#1}}

\newcommand{\QED}{\hfill $\Box$}

\def\P{\operatorname*{P}}

\def\Cov{\operatorname*{Cov}}
\def\E{\operatorname*{E}}
\def\I{\operatorname*{\mathbb{I}}}

\newcommand{\II}[1]{\mathbb{I}\left(#1\right)}

\newcommand{\nwc}{\newcommand}
\nwc{\COM}[1]{}

\def\IF{\infty}

\def\bs{\boldsymbol}

\begin{document}
\centerline{\bf \Large Berman's Inequality under Random Scaling}

\centerline{Enkelejd Hashorva and Zhichao Weng}

\centerline{University of Lausanne, Faculty of Business and Economics (HEC Lausanne),
Lausanne 1015, Switzerland}

\centerline{\today{}}

{\bf Abstract}: Berman's inequality is the key for establishing asymptotic properties of \tE{maxima of} Gaussian random sequences and supremum of Gaussian random fields.
This contribution shows that, asymptotically an extended version of Berman's inequality can be established for randomly scaled Gaussian random \Ke{vectors}. Two applications presented in this paper demonstrate the use of Berman's inequality under random scaling.

{\bf Key words}: Berman's inequality; Limit distribution; Extremal index;   Random scaling; H\"usler-Reiss distribution.

{\bf AMS Classification}: Primary 60G15; Secondary 60G70

\def\ARIJ{A_{ij}}
\def\DUV{\Delta_{\vk{S}}(\vk{u},\vk{v})}
\def\DUVs{\Delta_{S \vk{1}}(\vk{u},\vk{v})}

\section{Introduction}
\COM{\aH{An important contribution in extreme value theory concerned with maxima of triangular arrays of Gaussian random variables is \cite{MR1398063}. Motivated by the findings of
H\"usler and Reiss in 1989  (see \cite{hue1989}) the aforementioned contribution considered
a triangular array of \tE{$N(0,1)$} random variables $\{X_{n,i}, i,n \ge \wz{1}\}$ such that
for each $n$, $\{X_{n,i},i\ge 1\}$ is \tE{a stationary} Gaussian \Ke{random} sequence. Assuming that
$\wz{\varrho}_{n,j}=\EE{X_{n,i}X_{n,i+j}}$ satisfies some asymptotic conditions (see Theorem \ref{th3.1} below), \cite{MR1398063} shows that
\BQN\label{THR}
\lim_{n\to\IF}\pk{\max_{\wz{1}\le i \le n}X_{n,i}\le a_nx+ b_n}=\exp(-\vartheta\exp(-x)), \quad x\inr,
\EQN
where $\vartheta\in (0,1]$ is some known index and
\BQN\label{anbn}
 a_n = (2\ln n)^{-\frac{1}{2}} ,\quad b_n=(2\ln n)^{\frac{1}{2}}-\frac{1}{2}(2\ln n)^{-\frac{1}{2}}(\ln\ln n+\ln 4\pi).
 \EQN
The proof of \eqref{THR} strongly relies on Berman's inequality which is the basic
tool for the analysis of extreme values of Gaussian processes and \AH{Gaussian} random fields.}}
\zc{In the analysis of extreme values of Gaussian processes and \AH{Gaussian} random fields Berman's inequality plays a crucial role.}
 Essentially, for given two Gaussian distribution functions in $\R^d$ it bounds their difference by comparing the covariances. The key result that motivated \Ke{this} comparison method is  Plackett's partial differential equation given in \cite{Plackett}. As explained in \cite{Kratz06},
it was then developed by Slepian \cite{Slepian62}, Berman \cite{Berman64, Berman92}, Cram\'{e}r \cite{Cramer67}, Piterbarg \cite{Pit72, Pit96}
and then by Li and Shao \cite{LiShao02}. Specifically, the developed results are summarised by Berman's inequality which we formulate below in the \rH{most general form
derived} in \cite{LiShao02}. Let therefore \rH{$\vk{X}= (X_1 \ldot X_n)$ and $\vk{Y}=(Y_1 \ldot Y_n)$} be \rH{two Gaussian random vectors with $N(0,1)$ components} and covariance matrices $\cW{\Lambda_1}=(\lambda^{(1)}_{ij})$ and
$\cW{\Lambda_2}=(\lambda^{(2)}_{ij})$, respectively.
For arbitrary constants  $u_i, i\le n$, \cite{LiShao02} obtained
\BQNY
\pk{X_i\leq u_i,1\le i\le n }-\pk{Y_i \leq u_i, 1\le i\le n }
\leq \frac{1}{2 \pi }\sum_{1\leq i < j\leq n}
A_{ij}
\EXP{-\frac{u^2_{i}+u^2_{j}}{2(1+\rho_{ij})}},
\EQNY
where
\BQN\label{ARS}
\rho_{ij}:=\max(|\lambda^{(1)}_{ij}|,|\lambda^{(2)}_{ij}|), \quad A_{ij}=\big|\arcsin(\lambda^{(1)}_{ij})- \arcsin( \lambda^{(2)}_{ij})\big|.
\EQN
Clearly, for arbitrary constants $v_i,u_i, i\le n$, \aH{set $w:=\min_{1 \le i \le n} \min(\abs{u_i},\abs{v_i})$}
\BQN\label{ARIJ}
\pk{-v_i<X_i\leq u_i,1\le i\le n }-\pk{-v_i<Y_i \leq u_i, 1\le i\le n }
\leq \frac{2}{\pi }\sum_{1\leq i < j\leq n}
A_{ij}
\EXP{-\frac{w^2}{1+\rho_{ij}}},
\EQN
\COM{
where
\BQN
\ARIJ:=\zcw{\big|\arcsin(\lambda^{(1)}_{ij})- \arcsin( \lambda^{(2)}_{ij})\big|},
\EQN}
for a detailed proof see \cite{leadbetter1983extremes}, see also \cite{LuW2014} for recent extensions.\\ 
\cE{Berman's inequality can be applied also to non-Gaussian random \AH{vectors}. For instance, consider
two random vectors $\widetilde{\vk{X}}=(S_1 X_1 \ldot S_n X_n)$
and $\widetilde{\vk{Y}}=(S_1 Y_1 \ldot S_n Y_n)$ with
$S$, \rH{$S_{i}, i\le n$ some} positive \wz{independent} random variables with \cH{common} distribution function $G$ being further independent from \rH{$\vk{X}$ and $\vk{Y}$.}} \aH{In the special case}  $G$ is the uniform distribution on $(0,1)$, \aH{the upper bound in \eqref{ARIJ} implies}
\BQN\label{addeg1}
\DUV&:=& \cW{\pk{-v_i <S_iX_i\leq u_i, 1 \le i \le n }-\pk{-v_i<S_iY_i \leq u_i, 1 \le i \le n }}\nonumber\\
&\leq&\frac{2}{\pi}\sum_{1\leq i < j\leq n}
\ARIJ
\int_0^1\int_0^1\EXP{-\frac{(w/s_i)^2+(w/s_j)^2}{2(1+\rho_{ij})}}ds_i ds_j\nonumber\\
&\leq&\frac{2}{\pi}\sum_{1\leq i < j\leq n}\ARIJ\EXP{-\frac{w^2}{1+\rho_{ij}}}.
\EQN
 \aH{Another tractable case is when} $G(x)=1- \exp(-x),x>0$ is the exponential distribution. Indeed, by \eqref{ARIJ}
for all $0< a, b<1$
\aH{we have}
\BQN\label{addeg2}
\DUV&\leq&\frac{2}{\pi}\sum_{1\leq i < j\leq n}
\ARIJ
\int_0^{\infty}\int_0^{\infty}
\EXP{-\frac{(w/s_i)^2+(w/s_j)^2}{2(1+\rho_{ij})}-s_i-s_j}ds_ids_j\nonumber\\
&=&\frac{2}{\pi}\sum_{1\leq i < j\leq n}
\ARIJ
\int_0^{\infty}\int_0^{\infty}
\EXP{-\frac{(w/s_i)^2+(w/s_j)^2}{2(1+\rho_{ij})}-as_i-bs_j}\EXP{-(1-a)s_i-(1-b)s_j}ds_ids_j\nonumber\\
&\le&\frac{2}{\pi(1-a)(1-b)}\sum_{1\leq i < j\leq n}
\ARIJ\EXP{-\frac{3}{2}(a^{\frac{2}{3}}+b^{\frac{2}{3}})(1+\rho_{ij})^{-\frac{1}{3}}w^{\frac{2}{3}}}.
\EQN

\cE{Clearly, if we do not know the distribution function of $S$ it is not possible to obtain an explicit upper bound for $\DUV$. Since for the analysis of extremes of Gaussian random sequences or processes  \tE{Berman's} inequality is applied for large values of the  $u_i$'s and \tw{$v_i$'s} (see e.g., \cite{Pit96}), in this paper we are concerned with the derivation of Berman's inequality for \tE{some} general scaling random variable $S$ \tE{and all $u_i$'s \qE{and} $v_i$'s sufficiently large}.} \He{We shall consider two particular cases for the random vector $\vk{S}=(S_1 \ldot S_n)$, namely it has independent components, and it is comonotonic with $\vk{S}=(S \ldot S)=: S \vk{1}$. From the \tE{proofs} it can be seen that the joint dependence of $(S_i, S_j)$ for any pair $(i,j)$ is crucial; our results can be \tE{in fact} extended for \tE{certain} tractable dependence models. \tE{We} shall deal for simplicity only with these two cases.}
\\
 \AH{Since random scaling is a natural phenomena \qE{related} to \qE{the} time-value of money in finance,
measurement errors in experimental data, or physical constrains, the extension of Berman's inequality for inflated/deflated Gaussian random vectors is or certain interest also for statistical applications.}\\
\aH{Of course, Berman's inequality alone is not enough for extending \cite{MR1398063} to randomly scaled Gaussian triangular arrays; some additional results (see \cite{HashorvaWengB,HashorvaWengA}) which show that} for some tractable tail assumptions on $S$ the scaled random vector $\widetilde{\vk{X}}$ behaves similarly to  $\vk{X}$ are also important.
\aH{Specifically, we shall
deal with two large classes of random scaling: a) $S$ is a bounded random variable \tE{with}
a \Ee{tractable} tail behaviour at the right  endpoint of its distribution function, including \tE{in particular} the case that its survival function is regularly varying, and
b) $S$ is a Weibull-type random variable.}\\
\Ee{In view of our findings, several known results for Gaussian random sequences and processes can be extended to the scaled Gaussian framework; we shall demonstrate this with two representative applications}.

\rH{Organisation of the rest of the paper:} \aH{Section 2 presents Berman's inequality for scaled Gaussian random vectors. In Section 3 we display two applications, while the proofs are relegated to Section 4.}

\section{Main Results}
\aH{We consider first the case that $S$ is non-negative with distribution function $G$ which has right endpoint equal 1. Intuitively, large values of $S$ do not influence significantly large values of the product say $SX$ if $X$ is \wz{a} Gaussian random variable \wz{\tE{being} independent of $S$}. \tE{It turns} out that the following asymptotic upper bound}
\BQN\label{Potter}
\pk{S> 1- \cW{1/u}} \le c_A u^{\rH{-\tau}}
\EQN
\aH{valid} for all $u$ large and some $c_A>0, \Ke{\tau \ge 0}$ is sufficient for the derivation of a useful upper bound for $\DUV$ defined above. \\
\tE{A} canonical example of such $S$ is the beta random variable, which is a special case \AH{of a power-tail random variables \tE{$S$}}, namely
\BQN\label{eq:SA} \pk{S>1-\rH{1/u}}=(1+o(1))cu^{-\tau}, \quad u\to {\IF}
\EQN
\tE{holds} for some $c>0, \Ke{\tau \ge 0}$. \aH{Hereafter we set $w=\min_{1 \le i \le n} \min(\abs{u_i},\abs{v_i})$ and write $\DUVs$ \Ke{instead of
$\DUV$} if $\vk{S}=(S \ldot S)$.} \zw{Further write $\Delta_{\vk{S}}(u\vk{1})$ and $\Delta_{S\vk{1}}(u\vk{1})$ instead of $\DUV$ if all components of $\vk{v}$ equal $-\IF$,
$\vk{u}=(u \ldot u)=:u\vk{1}$ and the covariance matrix $\Lambda_2$ of $\vk{Y}$ is \tE{identity} matrix.}

\BT\label{th2.1}
Let $\vk{X},\widetilde{\vk{X}},\vk{Y},\widetilde{\vk{Y}}, S, S_i,i\le n$ be as above. If \eqref{Potter} holds,
then for all $u_i, v_i, 1\leq i \leq n $ large and \cW{$\epsilon>0$} we have
\BQN\label{Ba}
\DUV &\leq&\AH{(\mathbb{K}_A+\epsilon)}w^{-4\rH{\tau}}\sum_{1\leq i < j\leq n}
\ARIJ
(1+\rho_{ij})^{2\rH{\tau}}
\EXP{-\frac{w^2}{1+\rho_{ij}}}
\EQN
\wzc{and
\BQN\label{eq2.4}
\DUVs&\leq&\AH{(\mathbb{K}_A^*+\epsilon)}w^{-2\rH{\tau}}\sum_{1\leq i < j\leq n}
\ARIJ
(1+\rho_{ij})^{\rH{\tau}}
\EXP{-\frac{w^2}{1+\rho_{ij}}},
\EQN
}
where $\AH{\mathbb{K}_A}=\frac{2}{\pi}\cW{c_A^2}(\Gamma(\rH{\tau}+1))^2$ and \wzc{$\AH{\mathbb{K}_A^*}=\frac{2^{1-\tau}}{\pi}\cW{c_A}\Gamma(\rH{\tau}+1)$}.
\ET

\BK\label{coadd1}
\wz{Under the conditions of Theorem \ref{th2.1}, for all \wzc{u large and some positive constants $\mathcal{Q}$} we have
\BQN
\Delta_{\vk{S}}(\zw{u\vk{1}})\leq \mathcal{Q}u^{-4\rH{\tau}}\sum_{\wz{1\le i< j\le n}}|\wz{\lambda^{(1)}_{ij}}|
\EXP{-\frac{u^2}{1+|\wz{\lambda^{(1)}_{ij}}|}}
\EQN
and
\wzc{
\BQN
\Delta_{S\vk{1}}(u\vk{1})\leq\mathcal{Q}u^{-2\rH{\tau}}\sum_{\wz{1\le i< j\le n}}|\wz{\lambda^{(1)}_{ij}}|
\EXP{-\frac{u^2}{1+|\wz{\lambda^{(1)}_{ij}}|}}.
\EQN
}
}
\EK

\COM{
\BK\label{co2.1}
Under the conditions of Theorem \ref{th2.1}, for all $u$ large and some positive constants $\mathcal{Q}$ we have
\EK
}
\aH{We \Ee{shall investigate} \Ke{below the} more difficult case that the scaling random variable $S$ has distribution function with an infinite right endpoint.}
\aH{Motivated by the example of the exponential distribution in Introduction, we shall assume that $S$ has tail behaviour similar to a Weibull distribution}.  Specifically, for given constants $\alpha\inr, c_B,L, p \in (0,\infty)$ suppose that
\BQN \label{eq:SB}
\pk{S> u}=(1+o(1)) c_B u^{\alpha} \exp(-L u^{p}), \quad u\to \infty
 \EQN
is valid. \aH{The class of distribution functions satisfying \tE{\eqref{eq:SB} is quite} large. More importantly, \AH{under \eqref{eq:SB}} $S X$ has also  a \AH{\wz{Weibull} tail behaviour} if $X$ is a $N(0,1)$ random variable independent of $S$, see e.g., \cite{HashorvaWengA}. We state next our second result for Weibull-type random scaling.}

\COM{
\rH{Shall we suppose instead of \eqref{eq:SB} the following
\BQN \label{eq:SBB}
\pk{S> x} \le c_B u^{\alpha} \exp(-L u^{p})
 \EQN
holds for all large $u$???}
}

\BT\label{th2.2}
Let $\vk{X},\widetilde{\vk{X}},\vk{Y},\widetilde{\vk{Y}}, S, S_i,i\le n$ be as above. If \eqref{eq:SB} holds,
then for all $u_i, v_i, 1\leq i \leq n $ large and \cW{$\epsilon>0$}
 we have
\BQN \label{bb}
\DUV &\leq&
(\mathbb{K}_B+ \epsilon)w^{\frac{4\alpha+2p}{2+p}}\sum_{1\leq i < j\leq n}
\ARIJ
(1+\rho_{ij})^{\frac{-2\alpha-p}{p+2}}
\EXP{-2(1+\rho_{ij})^{-\frac{p}{2+p}}Tw^{\frac{2p}{2+p}}}
\EQN
\wzc{and
\BQN
\DUVs&\leq&(\mathbb{K}_B^*+ \epsilon)w^{\frac{2\alpha+p}{2+p}}\sum_{1\leq i < j\leq n}
\ARIJ
(1+\rho_{ij})^{\frac{-2\alpha-p}{2(p+2)}}
\EXP{-(2(1+\rho_{ij})^{-1})^{\frac{p}{2+p}}Tw^{\frac{2p}{2+p}}},
\EQN
}
where $T=L^{\frac{2}{p+2}}p^{-\frac{p}{p+2}}+(Lp)^{\frac{2}{p+2}}2^{-1}$, $\AH{\mathbb{K}_B}=4 \cW{c_B^2}(Lp)^{\frac{2(1-\alpha)}{p+2}}(p+2)^{-1}$ and
$\AH{\mathbb{K}_B^*}=2^{\frac{3+2p+\alpha}{2+p}}\pi^{-\frac{1}{2}}\cW{c_B}(Lp)^{\frac{1-\alpha}{p+2}}(p+2)^{-\frac{1}{2}}$.
\ET

\BK\label{coadd2}
\wz{Under the conditions of Theorem \ref{th2.2}, for all \wzc{$u$ large and some positive constants $\mathcal{Q}$} we have}
\BQN
\Delta_{\vk{S}}(u\vk{1})\leq \mathcal{Q}u^{\frac{4\alpha+2p}{2+p}}\sum_{\wz{1\le i <j\le n}}|\wz{\lambda^{(1)}_{ij}}|
\EXP{-2(1+|\wz{\lambda^{(1)}_{ij}}|)^{-\frac{p}{2+p}}T u^{\frac{2p}{2+p}}}
\EQN
and
\BQN
\Delta_{S\vk{1}}(u\vk{1})\leq \mathcal{Q}u^{\frac{2\alpha+p}{2+p}}\sum_{\wz{1\le i <j\le n}}|\wz{\lambda^{(1)}_{ij}}|
\EXP{-(2(1+|\wz{\lambda^{(1)}_{ij}}|)^{-1})^{\frac{p}{2+p}}Tu^{\frac{2p}{2+p}}}.
\EQN
\EK

\Ke{
{\bf Remarks}: a) Clearly, when $S$ is uniformly distributed on $(0,1)$ then condition (2.1) holds with $c_A=\tau=1$. For this case we have two results, the one derived in the Introduction and that given by \eqref{Ba}. We see that \zw{the bound obtained by \eqref{Ba} with $c_A=\tau=1$ is better due to the term $w^{-4\tau}$.}\\
b) Also for the case $S$ is a unit exponential random variables we have two bounds, one which holds for all values of $u_i,v_i,i\le n$ and the asymptotic one given in Theorem 2.3.
\zw{The bound implied by \eqref{bb} with $c_B=1, \alpha=0, p=1, L=1$ is asymptotically better than \tE{that implied by} \eqref{addeg2}.}
}

\COM{
\BK\label{co2.2}
Under the conditions of Theorem \ref{th2.2}, for all $u$ large and some positive constants $\mathcal{Q}$ we have
\BQNY
\left|\pk{S_iX_i\leq u, \mbox{for } i=1,\cdots,n}-\prod_{i=1}^n\pk{S_iX_i\leq u}\right|
\leq \mathcal{Q}u^{\frac{4\alpha+2p}{2+p}}\sum_{\wz{1\le i <j\le n}}|\wz{\lambda^{(1)}_{ij}}|
\EXP{-2(1+|\wz{\lambda^{(1)}_{ij}}|)^{-\frac{p}{2+p}}T u^{\frac{2p}{2+p}}}.
\EQNY
\EK
}

\COM{
{\bf Model A}: Consider the scaling random variable $S\ge 0$ with distribution function $G$ which has right endpoint equal 1.
Assume that there exists some random variable $
S_{\tau}$ which has a regularly varying survival function at 1 with index 
$\tau \in [0,\infty)$ such that
for any $u \in (\nu, 1)$ with $\nu \in (0,1)$  
\BQN\label{eq:bound S}
\pk{S_{\tau}> u}\geq \pk{S> u}.
\EQN
Suppose further that
$S_{\tau}$ have a regularly varying survival function at 1 with index 
${\tau}$. In view of Potter bound, there exist some $c_A$ and $\epsilon \in (0, \tau)$ such that
$$\pk{S_{\tau}> 1-x} \le c_A x^{\tau-\epsilon}.$$

If $S$ is a Beta random variable, then clearly $S$ satisfies  assumptions of Model A. 
Another instance is $S$ is such that $\pk{S=1}=c\in (0,1)$ and for some $\lambda <1$ we have $\pk{S<\lambda}=1- c$.
A more general example is if $S$ is such that 
that
\BQN\label{eq:SA} \pk{S>1-x}=(1+o(1))cx^{\gamma}, \quad x\to 0
\EQN
for some $c>0$.

{\bf Model B}: Assume that the  scaling random variable $S$ has an infinite right  endpoint, and its tail behaviour is similar to Weibull distributions.  More precisely,
for given constants $\alpha\inr, c_B,L, p \in (0,\infty)$ assume that
\BQN \label{eq:SB}
\pk{S> x}=(1+o(1)) c_B x^{\alpha} \exp(-L x^{p}), \quad x\to \infty
 \EQN
is valid.

\BT\label{th:SNCL m,M}
Suppose $\{X_i, 1\leq i \leq n\}$ are standard
Gaussian random variables with covariance matrix $\Lambda^{(1)}=(\lambda^{(1)}_{ij})$, and $\{Y_i, 1\leq i \leq n\}$
similarly with covariance matrix $\Lambda^{(2)}=(\lambda^{(2)}_{ij})$, and let $\rho_{ij}=\max(|\lambda^{(1)}_{ij}|,|\lambda^{(2)}_{ij}|)$.
Let $\{S_i, 1\leq i \leq n\}$ be iid random variables, and independent of $\{X_i\}$ and $\{Y_i\}$ respectively.
Further all $u_i, v_i, 1\leq i \leq n $ large and write
$w=\min(\mid u_1\mid, \ldots,\mid u_n\mid, \mid v_1\mid, \ldots,\mid v_n\mid)$, then
\begin{itemize}
\item[(1).] If \eqref{eq:bound S} holds we have 
\BQNY
&&\pk{-v_i<S_iX_i\leq u_i, \mbox{for } i=1,\cdots,n}-\pk{-v_i<S_iY_i \leq u_i, \mbox{for } i=1,\cdots,n}\\
&\leq&\mathbb{K}_A\sum_{1\leq i < j\leq n}
\ARIJ
(1+\rho_{ij})^{2\tau-2\epsilon}
w^{-4\tau+4\epsilon}
\EXP{-\frac{w^2}{1+\rho_{ij}}},
\EQNY
where $\mathbb{K}_A=\frac{2}{\pi}(c_A \Gamma(\tau-\epsilon+1))^2$.
\item[(2).] If \eqref{eq:SB} holds, then
\BQNY
&&\pk{-v_i<S_iX_i\leq u_i, \mbox{for } i=1,\cdots,n}-\pk{-v_i<S_iY_i \leq u_i, \mbox{for } i=1,\cdots,n}\\
&\leq&\mathbb{K}_B\sum_{1\leq i < j\leq n}
\ARIJ
(1+\rho_{ij})^{\frac{-2\alpha-p}{p+2}}
w^{\frac{4\alpha+2p}{2+p}}
\EXP{-2(1+\rho_{ij})^{-\frac{p}{2+p}}Tw^{\frac{2p}{2+p}}}
\EQNY
where $\mathbb{K}_B=4 c_B^2(Lp)^{\frac{2(1-\alpha)}{p+2}}(p+2)^{-1}$ and $T=L^{\frac{2}{p+2}}p^{-\frac{p}{p+2}}+(Lp)^{\frac{2}{p+2}}2^{-1}$.
\end{itemize}
\ET

\COM{
\BT\label{th:SNCL}
Suppose $\{X_i, 1\leq i \leq n\}$ are standard
normal variables with covariance matrix $\Lambda^{(1)}=(\lambda^{(1)}_{ij})$, and $\{Y_i, 1\leq i \leq n\}$
similarly with covariance matrix $\Lambda^{(2)}=(\lambda^{(2)}_{ij})$, and let $\rho_{ij}=\max(|\lambda^{(1)}_{ij}|,|\lambda^{(2)}_{ij}|)$.
Let $\{S_i, 1\leq i \leq n\}$ be iid random variables, and independent of $\{X_i\}$ and $\{Y_i\}$ respectively.
Further all $u_i$, $1\leq i\leq n$ large.
\begin{itemize}
\item[(1).] If \eqref{eq:bound S} holds, for some positive constants $c, \epsilon$, we have
\BQNY
&&\pk{S_iX_i\leq u_i, \mbox{for } i=1,\cdots,n}-\pk{S_iY_i \leq u_i, \mbox{for } i=1,\cdots,n}\\
&\leq&\mathcal{Q}_A\sum_{1\leq i < j\leq n}
\ARIJ
(1+\rho_{ij})^{2\tau-2\epsilon}
(u_{i}u_{j})^{-2\tau+2\epsilon}
\EXP{-\frac{u^2_{i}+u^2_{j}}{2(1+\rho_{ij})}},
\EQNY
where $\mathcal{Q}_A=\frac{1}{2\pi}(c\Gamma(\tau-\epsilon+1))^2$.
\item[(2).] If \eqref{eq:SB} holds we have
\BQNY
&&\pk{S_iX_i\leq u_i, \mbox{for } i=1,\cdots,n}-\pk{S_iY_i \leq u_i, \mbox{for } i=1,\cdots,n}\\
&\leq&\mathcal{Q}_B\sum_{1\leq i < j\leq n}
\ARIJ
(1+\rho_{ij})^{\frac{-2\alpha-p}{p+2}}
(u_{i}u_{j})^{\frac{2\alpha+p}{2+p}}
\EXP{-(1+\rho_{ij})^{-\frac{p}{2+p}}T\left(u_i^{\frac{2p}{2+p}}+u_j^{\frac{2p}{2+p}}\right)}
\EQNY
where $\mathcal{Q}_B=C^2(Lp)^{\frac{2(1-\alpha)}{p+2}}(p+2)^{-1}$ and $T=L^{\frac{2}{p+2}}p^{-\frac{p}{p+2}}+(Lp)^{\frac{2}{p+2}}2^{-1}$.
\end{itemize}
\ET
}

\BK\label{co:SNCL}
Assume that $\{X_i,1\leq i\leq n\}$ is a centered stationary Gaussian random sequence with
$\lambda_i=\Cov(X_1,X_{1+i})$.
Let $\{S_i, 1\leq i \leq n\}$ be iid random variables, and independent of $\{X_i\}$. Then for all $u$ large and some positive constants $\mathcal{Q}$
\begin{itemize}
\item[(1).] If \eqref{eq:bound S} holds, then
\BQNY
\left|\pk{S_iX_i\leq u, \mbox{for } i=1,\cdots,n}-\prod_{i=1}^n\pk{S_iX_i\leq u}\right|
\leq \mathcal{Q}n\sum_{i=1}^{n}|\lambda_i|
u^{-4\tau+4\epsilon}
\EXP{-\frac{u^2}{1+|\lambda_i|}}.
\EQNY
\item[(2).] Assuming \eqref{eq:SB}  we obtain
\BQNY
\left|\pk{S_iX_i\leq u, \mbox{for } i=1,\cdots,n}-\prod_{i=1}^n\pk{S_iX_i\leq u}\right|
\leq \mathcal{Q}n\sum_{i=1}^{n}|\lambda_i|
u^{\frac{4\alpha+2p}{2+p}}
\EXP{-2(1+|\lambda_i|)^{-\frac{p}{2+p}}T u^{\frac{2p}{2+p}}}
\EQNY
where $T=L^{\frac{2}{p+2}}p^{-\frac{p}{p+2}}+(Lp)^{\frac{2}{p+2}}2^{-1}$.
\end{itemize}
\EK

\textbf{Example 1.} Suppose $\{X_i, 1\leq i \leq n\}$ and $\{Y_i, 1\leq i \leq n\}$ defined as Theorem \ref{th:SNCL m,M}.
Let $\{S_i, 1\leq i \leq n\}$ be iid random variables with uniformly distribution, i.e., $\pk{S_1>1-x}=x$,
and independent of $\{X_i\}$ and $\{Y_i\}$ respectively. Then for all large $u_i, v_i, 1\leq i\leq n$,
\BQNY
&&|\pk{-v_i <S_iX_i\leq u_i, \mbox{for } i=1,\cdots,n}-\pk{-v_i<S_iY_i \leq u_i, \mbox{for } i=1,\cdots,n}|\\
&\leq&\frac{2}{\pi}\sum_{1\leq i < j\leq n}
\ARIJ
\int_0^1\int_0^1\EXP{-\frac{(w/s_i)^2+(w/s_j)^2}{2(1+\rho_{ij})}}ds_ids_j\\
&\leq&\frac{2}{\pi}\sum_{1\leq i < j\leq n}
\ARIJ
(1+\rho_{ij})^{2}
w^{-4}
\EXP{-\frac{w^2}{1+\rho_{ij}}}.
\EQNY

\textbf{Example 2.} Suppose $\{X_i, 1\leq i \leq n\}$ and $\{Y_i, 1\leq i \leq n\}$ defined as Theorem \ref{th:SNCL m,M}.
Let $\{S_i, 1\leq i \leq n\}$ be iid random variables with exponential distribution, i.e., $\pk{S_1>x}=e^{-x}$,
and independent of $\{X_i\}$ and $\{Y_i\}$ respectively. Then for all large $u_i,v_i, 1\leq i\leq n$,
\BQNY
&&|\pk{-v_i<S_iX_i\leq u_i, \mbox{for } i=1,\cdots,n}-\pk{-v_i<S_iY_i \leq u_i, \mbox{for } i=1,\cdots,n}|\\
&\leq&\frac{2}{\pi}\sum_{1\leq i < j\leq n}
\ARIJ
\int_0^{\infty}\int_0^{\infty}
\EXP{-\frac{(w/s_i)^2+(w/s_j)^2}{2(1+\rho_{ij})}-s_i-s_j}ds_ids_j\\
&\leq&\frac{4}{3}\sum_{1\leq i < j\leq n}
\ARIJ
(1+\rho_{ij})^{-\frac{1}{3}}
w^{\frac{2}{3}}
\EXP{-3(1+\rho_{ij})^{-\frac{1}{3}}w^{\frac{2}{3}}}.
\EQNY

}

\section{Applications}
An important contribution in extreme value theory concerned with maxima of triangular arrays of Gaussian random variables is \cite{MR1398063}. Motivated by the findings of
H\"usler and Reiss in 1989  (see \cite{hue1989}) the aforementioned contribution considered
a triangular array of \tE{$N(0,1)$} random variables $\{X_{n,i}, i,n \ge \wz{1}\}$ such that
for each $n$, $\{X_{n,i},i\ge 1\}$ is \tE{a stationary} Gaussian \Ke{random} sequence. Assume that
$\wz{\varrho}_{n,j}=\EE{X_{n,i}X_{n,i+j}}$
\aH{satisfies} \qE{for any $j\ge 1$}
\begin{eqnarray}\label{eq1}
\lim_{n\to \IF}(1-\wz{\varrho}_{n,j})\ln n \qE{=}\ \delta_j \in (0,\IF\zw{)}, \quad \delta_0:=0 
\end{eqnarray}
and for each $n$, $\wz{\varrho}_{n,j}$ decays fast enough as $j$ increases. \qE{Under some additional conditions (see Theorem \ref{th3.1} below)
the deep contribution} 
 \cite{MR1398063} shows that \qE{for the maxima $M_n= \max_{\wz{1}\le i \le n}X_{n,i}$}
\BQN\label{THR}
\lim_{n\to\IF}\pk{M_n \le a_nx+ b_n}=\exp(-\vartheta\exp(-x)), \quad x\inr,
\EQN
where \BQN\label{anbn}
 a_n = (2\ln n)^{-\frac{1}{2}} ,\quad b_n=(2\ln n)^{\frac{1}{2}}-\frac{1}{2}(2\ln n)^{-\frac{1}{2}}(\ln\ln n+\ln 4\pi)
 \EQN
and 
$$\vartheta=\pk{E/2+\sqrt{\delta_{\zw{k-1}}}W_k\le \delta_{\zw{k-1}}, \ \  \mbox{for all}\ \  k\ge \zw{2} },$$
with $E$ a unit exponential random variable independent of
$W_k$ and $\{W_k, k\ge \zw{2}\}$ being jointly Gaussian  with zero means and covariances
$$\EE{W_iW_j}=\frac{\delta_{\zw{i-1}}+\delta_{\zw{j-1}}-\delta_{|i-j|}}{2\sqrt{\delta_{i-1}\delta_{j-1}}}. 
$$
The proof of \eqref{THR} strongly relies on Berman's inequality. 
\zc{Hence,}
our first application \aH{extends} the result of \cite{MR1398063} to \aH{triangular arrays} of \tE{randomly} scaled Gaussian random variables.
\AH{In the following we \Ke{investigate} the effect of a comonotonic random scaling considering a bounded $S$ and thus \Ke{$\vk{S}=S \vk{1}$}.}

\begin{theorem}\label{th3.1}
Let $\{X_{n,i},i,n\ge 1\}$ be \zc{a triangular array of standard Gaussian random variables defined as above satisfying \eqref{eq1}}, \AH{being further independent of  the iid non-negative random variables $\{S_n,n\ge 1\}$ where $S_1$ satisfies \eqref{eq:SA}.} If there exist positive integers $r_n,l_n$ such that
\begin{eqnarray}\label{add3.1}
\lim_{n\to \IF}\frac{l_n}{r_n}=0,\quad \lim_{n\to\IF}\frac{r_n}{n}= 0,
\end{eqnarray}
\begin{eqnarray}\label{add3.2}
\lim_{n\to\IF}\frac{n^2}{r_n}c_n^{-\zw{\tau}}\sum_{j=l_n}^n  \frac{|\wz{\varrho}_{n,j}| (1+\zw{|\varrho_{n,j}|})^{\wz{\tau}}}
{\sqrt{ 1-\wz{\varrho}_{n,j}^2}} \exp\left(-\frac{c_n}{1+\cW{|\wz{\varrho}_{n,j}|}}\right)=0, \quad c_n:=2\ln n-(2\tau+1)\ln\ln n
\end{eqnarray}
and further
\begin{equation}\label{add3.3}
\lim_{m\to \IF}\limsup_{n\to \IF}\sum_{j=m}^{r_n}n^{-\frac{1-\wz{\varrho}_{n,j}}{1+\wz{\varrho}_{n,j}}}\frac{(\ln n)^{\frac{\tau(1-\wz{\varrho}_{n,j})-\wz{\varrho}_{n,j}}{1+\wz{\varrho}_{n,j}}}}{\sqrt{1-\wz{\varrho}_{n,j}^2}}=0,
\end{equation}
\qE{then for the maxima $M_n=\max_{\cW{1}\le i \le n}S_nX_{n,i} $ the result in \eqref{THR} holds with $\vartheta$ defined as above and}
\begin{equation}\label{anbntau}
a_n=(2\ln n)^{-1/2},\quad  \quad b_n=(2\ln n)^{1/2}+(2\ln n)^{-1/2}\left(\ln (c(2\pi)^{-1/2}\Gamma(1+\tau))-\frac{2\tau+1}{2}(\ln\ln n+\ln2)\right).
\end{equation}
\end{theorem}

\zc{{\bf Remark}: Using similar \qE{arguments} as \qE{in the proof of } Theorem \ref{th3.1}, the findings of the recent contribution \cite{DavisFrench} can also be extended by considering randomly scaled Gaussian field on a lattice.}

In our second application we consider scaled Gaussian random vectors where the scaling vector $\vk{S}$ has
independent components.
Let $\bigl\{\mathbf{X}_{n,k}=\bigl(X_{n,k}^{(1)},X_{n,k}^{(2)}\bigr), 1\leq k\leq n, n\geq 1\bigr\}$ \tE{be}
a triangular array of bivariate centered \wz{stationary} Gaussian random vectors with unit-variance and correlation given by
\begin{eqnarray*}
corr\left(X_{n,k}^{(1)},X_{n,k}^{(2)}\right)=\lambda_0(n),   \qquad
corr\left(X_{n,k}^{(i)},X_{n,l}^{(j)}\right)
\wz{=\lambda_{ij}(|k-l|,n)},
\end{eqnarray*}
where $1\leq k\neq l\leq n$ and $i,j \in \{1,2\}.$
Further, let $\bigl\{\hat{\mathbf{X}}_{n,k}=\bigl(\hat{X}_{n,k}^{(1)},\hat{X}_{n,k}^{(2)}\bigr), 1\leq k\leq n, n\geq 1\bigr\}$
denote the associated iid triangular array of $\{\mathbf{X}_{n,k}\}$, i.e.,
$corr\left(\hat{X}_{n,k}^{(1)},\hat{X}_{n,k}^{(2)}\right)=\lambda_0(n)$ and
$corr\left(\hat{X}_{n,k}^{(i)},\hat{X}_{n,l}^{(j)}\right)=0,$ for $1\leq k\neq l\leq n$ and $i,j \in \{1,2\}.$
 Let $\{S_{n,k}, 1\le k \le n, n\geq 1\}$ be iid
random variables \tE{being} independent of
$\{\mathbf{X}_{n,k}, 1\le k\le n, n\ge 1\}$ and $\{\hat{\mathbf{X}}_{n,k},1\le k\le n, n\ge 1\}$, respectively.
\COM{
For notational simplicity, in the sequel we shall use the following notation:
$$\left(M_{n}^{(1)},M_{n}^{(2)}\right):=\left(\max_{1\leq k\leq n}S_{nk}X_{n,k}^{(1)},\max_{1\leq k\leq n}S_{nk}X_{n,k}^{(2)}\right),
\qquad \left(m_{n}^{(1)},m_{n}^{(2)}\right):=\left(\min_{1\leq k\leq n}S_{nk}X_{n,k}^{(1)},\min_{1\leq k\leq n}S_{nk}X_{n,k}^{(2)}\right),$$
$$\left(\mathfrak{M}_{n}^{(1)},\mathfrak{M}_{n}^{(2)}\right):=\left(\max_{1\leq k\leq n}S_{nk}\hat{X}_{nk}^{(1)},\max_{1\leq k\leq n}S_{nk}\hat{X}_{nk}^{(2)}\right),
\qquad \left(\mathfrak{m}_{n}^{(1)},\mathfrak{m}_{n}^{(2)}\right):=\left(\min_{1\leq k\leq n}S_{nk}\hat{X}_{nk}^{(1)},\min_{1\leq k\leq n}S_{nk}\hat{X}_{nk}^{(2)}\right).$$
}
Assume that the correlation $\lambda_0(n)$
satisfies
\begin{eqnarray}\label{add 1}
\lim_{n\to \infty}\frac{b_n}{a_n}(1-\lambda_{0}(n)) = 2\lambda^2 \quad \mbox{with} \ \  \lambda\in [0,\infty],
\end{eqnarray}
where
$$\Ke{a_n=\frac{1}{1-F(b_n)} \int_{b_n}^{\Ee{\IF}}(1-F(x))dx}, \quad
b_n=F^{-1}\left(1-\frac{1}{n}\right),$$
with $F^{-1}$ the inverse of the df \Ke{$F$}
of $S_{1,1}\hat{X}_{1,1}^{(1)}$.
 \Ke{It is well-known (see e.g.,  \cite{HB}) that}
\BQNY\label{eq1.4}
\lim_{n \to
\infty}\sup_{x,y \in \R}\left|\pk{\max_{1\le k \le n}S_{n,k}\hat{X}_{n,k}^{(1)}\leq
u_n(x),\max_{1\le k \le n}S_{n,k}\hat{X}_{n,k}^{(2)}\leq u_n(y)}-H_{\lambda}(x,y)\right|=0,
\EQNY
where $u_n(z)=a_nz+b_n, z \in \R$ \AH{and the H\"usler-Reiss distribution function $H_\lambda$ is given by}
 \BQN \label{add 2}
H_{\lambda}(x,y)=\exp\left(-
e^{-x} \Phi\left(\lambda+\frac{y-x}{2\lambda}\right)-e^{-y} \Phi\left(\lambda+\frac{x-y}{2\lambda}\right)\right),
\EQN
with $\Phi$ the distribution function of an $N(0,1)$ random variable.

\AH{In the following we are interested in the case that only a fraction of random vectors is observed. Assume therefore that
 $\varepsilon_{n,k}$ is an indicator random variable of the event that \tE{the random vector} $\mathbf{X}_{n,k}$
is observed. Then $\Xi_n=\sum_{k=1}^n \varepsilon_{n,k}$ is the number of observed random vectors from the
set $\{\mathbf{X}_{n,1}, \cdots, \mathbf{X}_{n,n}\}$.\\
We shall impose the following condition:}

{\bf Condition E}. \AH{The indicator random variables $\varepsilon_{n,k}$ are independent of $\mathbf{X}_{n,k}$ and $S_{n,k}$.
Further, the convergence in probability}
$$\frac{\Xi_{n}}{n} \stackrel{P}{\to} \eta , \quad n\to \IF$$
\AH{holds with  $\eta$ some random variable taking values in $(0,1]$ almost surely.}

\wz{For \Ee{notational simplicity we set}
\begin{eqnarray*}
\mathbf{M}_n(\bs{\varepsilon}_n):=
\left\{
\begin{array}{ll}
\max\{S_{n,k}\mathbf{X}_{n,k}, 1\le k\le n, \varepsilon_{n,k}=1 \}, & if \  \sum_{k=1}^n\varepsilon_{n,k}\ge 1,\\
\inf\{\mathbf{x}|\pk{S_{n,k}\mathbf{X}_{n,k}\leq\mathbf{x}}>\mathbf{0}\}, & otherwise,\\
\end{array}
\right.
\end{eqnarray*}
\begin{eqnarray*}
\mathbf{m}_n(\bs{\varepsilon}_n):=
\left\{
\begin{array}{ll}
\min\{S_{n,k}\mathbf{X}_{n,k}, 1\le k\le n, \varepsilon_{n,k}=1 \}, & if \ \sum_{k=1}^n\varepsilon_{n,k}\ge 1,\\
\inf\{\mathbf{x}|\pk{S_{n,k}\mathbf{X}_{n,k}\leq\mathbf{x}}>\mathbf{0}\}, & otherwise\\
\end{array}
\right.
\end{eqnarray*}
and $\mathbf{M}_n=\max\{S_{n,k}\mathbf{X}_{n,k},1\le k\le n\}$, $\mathbf{m}_n=\min\{S_{n,k}\mathbf{X}_{n,k},1\le k\le n\}.$}
\COM{For subset $I_n \subset N$, let us
define the following notation:
\begin{eqnarray*}
\mathbf{M}(I_n, \bs{\varepsilon}_n):=
\left\{
\begin{array}{ll}
\max\{S_{n,k}\mathbf{X}_{n,k}: k\in I_n, \varepsilon_{n,k}=1 \}, & if \  \sum_{k \in I_n}\varepsilon_{n,k}\ge 1;\\
\inf\{\mathbf{x}|\pk{S_{n,k}\mathbf{X}_{n,k}\leq\mathbf{x}}>\mathbf{0}\}, & otherwise.\\
\end{array}
\right.
\end{eqnarray*}
\begin{eqnarray*}
\mathbf{m}(I_n,\bs{\varepsilon}_n):=
\left\{
\begin{array}{ll}
\min\{S_{n,k}\mathbf{X}_{n,k}, k\in I_n, \varepsilon_{n,k}=1 \}, & if \ \sum_{k \in I_n}\varepsilon_{n,k}\ge 1;\\
\inf\{\mathbf{x}|\pk{S_{n,k}\mathbf{X}_{n,k}\leq\mathbf{x}}>\mathbf{0}\}, & otherwise.\\
\end{array}
\right.
\end{eqnarray*}
\begin{eqnarray*}
\widehat{\mathbf{M}}(I_n,\bs{\varepsilon}_n):=
\left\{
\begin{array}{ll}
\max\{S_{n,k}\hat{\mathbf{X}}_{n,k}: k \in I_n, \varepsilon_{n,k}=1 \}, & if \ \sum_{k \in I_n}\varepsilon_{n,k}\ge 1;\\
\inf\left\{\mathbf{x}\big{|}\pk{S_{n,k}\hat{\mathbf{X}}_{n,k}\leq\mathbf{x}}>\mathbf{0}\right\}, & otherwise.\\
\end{array}
\right.
\end{eqnarray*}
\begin{eqnarray*}
\widehat{\mathbf{m}}(I_n,\bs{\varepsilon}_n):=
\left\{
\begin{array}{ll}
\min\{S_{n,k}\hat{\mathbf{X}}_{n,k}, k\in I_n, \varepsilon_{n,k}=1 \}, & if \ \sum_{k \in I_n}\varepsilon_{n,k}\ge 1;\\
\inf\left\{\mathbf{x}\big{|}\pk{S_{n,k}\hat{\mathbf{X}}_{n,k}\leq\mathbf{x}}>\mathbf{0}\right\}, & otherwise.\\
\end{array}
\right.
\end{eqnarray*}
For simplicity, we write
$\mathbf{M}_n(\bs{\varepsilon}_n)=\mathbf{M}(\{1,2,\ldots,n\}, \bs{\varepsilon}_n),$
$\mathbf{M}(I_n)=\max\{S_{n,k}\mathbf{X}_{n,k}, k\in I_n\},$
$\mathbf{M}_n=\max\{S_{n,k}\mathbf{X}_{n,k},1\le k\le n\}.$
Similarly we also define $\mathbf{m}_n(\bs{\varepsilon}_n),
\mathbf{m}(I_n), \mathbf{m}_n, \widehat{\mathbf{M}}_n(\bs{\varepsilon}_n),
\widehat{\mathbf{M}}(I_n), \widehat{\mathbf{M}}_n, \widehat{\mathbf{m}}_n(\bs{\varepsilon}_n),
\widehat{\mathbf{m}}(I_n), \widehat{\mathbf{m}}_n.$}
\COM{
For random variable $\eta$ such that $0\le \eta \le 1$ a.s., write
\begin{eqnarray*}
B_{\mu,l}=
\left\{\omega: \eta(\omega) \in \left\{
\begin{array}{ll}
[0,\frac{1}{2^l}], & \mu=0,\\
(\frac{\mu}{2^l},\frac{\mu+1}{2^l}], & 0< \mu \le 2^l-1\\
\end{array}
\right\},\right.
\end{eqnarray*}
$$B_{\mu, l, \bs{\beta}_n}=\{\omega: \delta_{nk}(\omega)=\beta_{nk}, 1\le k\le n\}\cap B_{\mu, l}.$$
}
\COM{
Furthermore, define
\begin{array}{ll}
\mathbb{A}^{(1)}(w)=\{M_n^{(1)} \le u_n(w)\}, & \mathbb{A}^{(2)}(w)=\{M_n^{(1)} \le u_n(w)\} \\
\mathbb{A}^{(1)}(w)=\{M_n^{(1)} \le u_n(w)\}, & \mathbb{A}^{(2)}(w)=\{M_n^{(1)} \le u_n(w)\} \\
\mathbb{A}^{(1)}(w)=\{M_n^{(1)} \le u_n(w)\}, & \mathbb{A}^{(2)}(w)=\{M_n^{(1)} \le u_n(w)\} \\
\mathbb{A}^{(1)}(w)=\{M_n^{(1)} \le u_n(w)\}, & \mathbb{A}^{(2)}(w)=\{M_n^{(1)} \le u_n(w)\} \\
\mathbb{B}^{(1)}(w)=\{m_n^{(1)} > -u_n(w)\}, & \mathbb{B}^{(2)}(w)=\{m_n^{(1)} > -u_n(w)\} \\
\mathbb{B}^{(1)}(w)=\{m_n^{(1)} \le u_n(w)\}, & \mathbb{B}^{(2)}(w)=\{m_n^{(1)} \le u_n(w)\} \\
\mathbb{B}^{(1)}(w)=\{m_n^{(1)} \le u_n(w)\}, & \mathbb{B}^{(2)}(w)=\{m_n^{(1)} \le u_n(w)\} \\
\mathbb{B}^{(1)}(w)=\{m_n^{(1)} \le u_n(w)\}, & \mathbb{B}^{(2)}(w)=\{m_n^{(1)} \le u_n(w)\} \\
\end{array}
}

\AH{For $S_{n,k}=1, 1\leq k\leq n$ almost surely, 
\zw{according to \cite{HashorvaPengWengL}}, under Condition {\bf E} we have}
\BQNY
\lim_{n \to \infty}\sup_{x_1,y_1\in \R \atop x_1\le y_1}
\left|\pk{M^{(1)}_n(\bs{\varepsilon}_n)\leq u_n(x_1),M^{(1)}_n\leq u_n(y_1)}-\EE{\zw{\Lambda}^{\eta}(x_1)\zw{\Lambda}^{1-\eta}(y_1)}\right|=0,
\EQNY
where $u_n(x)=a_nx+b_n$ with $a_n$ and $b_n$ defined in \eqref{anbn} \zw{and $\Lambda(x)=\EXP{-e^{-x}}$, $x \inr$, provided that
$\lim_{n \to \IF}\max_{l_n <k <n }\lambda_{11}(k,n)\ln n =0$ with $l_n=[n^{\hat{\beta}}]$,
$0<\hat{\beta}<(1-\hat{\sigma})/(1+\hat{\sigma})$ and $\hat{\sigma}=\max_{1\le k<n,n\ge1}|\lambda_{11}(k,n)|$.}
\AH{Below we obtain a more general result for our 2-dimensional setup considering Weibull-type random scaling.}

\begin{theorem}\label{th:m,M}
 Let $\{(X^{(1)}_{n,k}, X^{(2)}_{n,k}), 1\leq k\leq n,
n\geq 1\}$ be a bivariate triangular array of standard Gaussian  random vectors defined as above.
Let $\{S_{n,k},1\leq k\leq n, n\geq 1\}$ be iid random variables \Ke{being} independent of
$\{(X_{n,k}^{(1)},X_{n,k}^{(2)}), 1\leq k\leq n,
n\geq 1\}$. \tE{Suppose that} the correlation $\lambda_0(n)$ satisfy \eqref{add 1} with $\lambda\in(0,\IF)$ \Ke{and condition {\bf E} holds}.
 Let $\beta$ be a constant satisfying $0<\beta<2(1+\sigma)^{-\frac{p}{2+p}}-1$
with $\sigma=\max_{1\leq k<n, n\ge 1 \atop i,j \in \{1,2\}}|\lambda_{ij}(k,n)|<1$, and write $\iota_n=[n^{\beta}]$. If
\eqref{eq:SB} holds and the covariance function satisfies
\BQNY\label{add 3}
\lim_{n \to \infty}\max_{\iota_n\leq k<n \atop i,j \in \{1,2\}}\lambda_{ij}(k,n)\ln n=0,
\EQNY
then we have
\BQNY \label{eq2.2}
\lim_{n \to \infty}\sup_{x_i,y_i \in \R, i\le 4 \atop x_1\le x_3, x_2\le x_4, \zw{y_1 \le y_3,y_2 \le y_4}}&&
\left|\P\left(-u_n(y_1)<m^{(1)}_n(\bs{\varepsilon}_n)\leq M^{(1)}_n(\bs{\varepsilon}_n)\leq
u_n(x_1),-u_n(y_2)<m^{(2)}_n(\bs{\varepsilon}_n)\leq M^{(2)}_n(\bs{\varepsilon}_n)\leq u_n(x_2),\right.\right. \nonumber\\
&&\quad \left.-u_n(y_3)<m_n^{(1)}\leq M_n^{(1)}\leq
u_n(x_3),-u_n(y_4)<m_n^{(2)}\leq M_n^{(2)}\leq u_n(x_4)\right) \nonumber\\
&&\quad \left.-\E\left(H_{\lambda}^{\eta}(x_1,x_2)H_{\lambda}^{\eta}(y_1,y_2)H_{\lambda}^{1-\eta}(x_3,x_4)H_{\lambda}^{1-\eta}(y_3,y_4)\right)\right|=0,
\EQNY
where
$H_{\lambda}$ \AH{is} defined in \eqref{add 2} and norming
constants
$a_n$ and $b_n$ satisfy
\BQN\label{anbn B}
a_n=\frac{2+p}{2p}T^{-\frac{2+p}{2p}}(\ln
n)^{\frac{2-p}{2p}}, \qquad b_n=\left(\frac{\ln
n}{T}\right)^{\frac{2+p}{2p}}
+\frac{2+p}{2p}T^{-\frac{2+p}{2p}}(\ln
n)^{\frac{2-p}{2p}}\left(\frac{\alpha}{p}\ln \ln n
-\frac{\alpha}{p}\ln T+\ln \varpi_B\right) \EQN
with
$T=2^{-1}Q^2+LQ^{-p}$, $\varpi_B=c_B(2+p)^{-\frac{1}{2}}Q^{-\alpha}$
and $Q=(pL)^{1/(2+p)}$.
\end{theorem}

\COM{
\begin{itemize}
\item[(1).]
Suppose that \eqref{Potter} holds, and there exists some $\gamma \in [0,\infty)$ such that for all $u$ large and some $\tilde{c}_A>0$
\BQN\label{eq S gamma}
\cW{\pk{S>1-1/u}\geq \tilde{c}_A u^{-\gamma}}.
\EQN
Further let $r$ be a constant
such that $0< r<(1-\sigma)/(1+\sigma)$ and write $\kappa_n=[n^{r}]$,  assume that there exist some small $\epsilon>0$ such that
\begin{eqnarray}\label{eq2.1}
\lim_{n \to \infty}\max_{\kappa_n\leq k<n \atop i,j \in \{1,2\}}\lambda_{ij}(k,n)(\ln n)^{1+2(\gamma-\tau)}=0.
\end{eqnarray}
\COM{
\cW{Suppose that \eqref{eq:SA} holds, \COM{ and there exists some $\gamma \in [0,\infty)$ such that for all $u$ large and some $\tilde{c}_A>0$
\BQN\label{eq S gamma}
\cW{\pk{S>1-1/u}\geq \tilde{c}_A u^{-\gamma}}.
\EQN
Further}
 let $r$ be a constant
such that $0< r<(1-\sigma)/(1+\sigma)$ and write $\kappa_n=[n^{r}]$,
and the covariance satisfy
\begin{eqnarray}\label{eq2.1}
\lim_{n \to \infty}\max_{\kappa_n\leq k<n \atop i,j \in \{1,2\}}\lambda_{ij}(k,n)\ln n=0,
\end{eqnarray}}}
Then we have
\BQN \label{eq2.2}
\lim_{n \to \infty}\sup_{x_i,y_i \in \R, i\le 4 \atop x_1\le x_3, x_2\le x_4, y_1 \ge y_3,y_2 \ge y_4}&&
\left|\P\left(-u_n(y_1)<m^{(1)}_n(\bs{\varepsilon}_n)\leq M^{(1)}_n(\bs{\varepsilon}_n)\leq
u_n(x_1),-u_n(y_2)<m^{(2)}_n(\bs{\varepsilon}_n)\leq M^{(2)}_n(\bs{\varepsilon}_n)\leq u_n(x_2),\right.\right. \nonumber\\
&&\quad \left.-u_n(y_3)<m_n^{(1)}\leq M_n^{(1)}\leq
u_n(x_3),-u_n(y_4)<m_n^{(2)}\leq M_n^{(2)}\leq u_n(x_4)\right) \nonumber\\
&&\quad \left.-\E\left(H_{\lambda}^{\eta}(x_1,x_2)H_{\lambda}^{\eta}(y_1,y_2)H_{\lambda}^{1-\eta}(x_3,x_4)H_{\lambda}^{1-\eta}(y_3,y_4)\right)\right|=0,
\EQN
where
$H_{\lambda}$ \AH{is} defined in \eqref{add 2} and norming
constants $a_{n}$ and $b_{n}$ satisfy \AH{(we can give exact formulas for $a_n, b_n$???}
\BQN\label{eq:ab}
\lim_{n \to \infty}b_n/\sqrt{2\ln n}=\lim_{n \to \infty}a_n\sqrt{2\ln n}=1.
\EQN
\item[(2).] Suppose that \eqref{eq:SB} holds, Set $\sigma=\max_{1\leq k<n \atop i,j \in \{1,2\}}|\lambda_{ij}(k,n)|<1$ and $\beta$ be a constant such that $0<\beta<2(1+\sigma)^{-\frac{p}{2+p}}-1$
and write $\iota_n=[n^{\beta}]$,
and the covariance satisfy
\BQN\label{add 3}
\lim_{n \to \infty}\max_{\iota_n\leq k<n \atop i,j \in \{1,2\}}\lambda_{ij}(k,n)\ln n=0,
\EQN
then \eqref{eq2.2}
also holds with $a_n$ and $b_n$ satisfying
\BQN\label{anbn B}
a_n=\frac{2+p}{2p}T^{-\frac{2+p}{2p}}(\ln
n)^{\frac{2-p}{2p}}, \qquad b_n=\left(\frac{\ln
n}{T}\right)^{\frac{2+p}{2p}}
+\frac{2+p}{2p}T^{-\frac{2+p}{2p}}(\ln
n)^{\frac{2-p}{2p}}\left(\frac{\alpha}{p}\ln \ln n
-\frac{\alpha}{p}\ln T+\ln \varpi_B\right), \EQN where
$T=2^{-1}Q^2+LQ^{-p}$, $\varpi_B=c_B(2+p)^{-\frac{1}{2}}Q^{-\alpha}$
and $Q=(pL)^{1/(2+p)}$. \end{itemize}
}

\COM{
\BT\label{th:in M}
 Under the conditions of Theorem \ref{th:m,M}, and assume that the conditions $C_1$ and $C_2$ hold.
\begin{itemize}
\item[(1).] Suppose that \eqref{eq:bound S}, \eqref{eq S gamma} and \eqref{eq2.1} hold,
then we have
\BQN \label{eq22}
\lim_{n \to \infty}\sup_{x_1,x_2,y_1,y_2 \in \R \atop x_1<y_1, x_2<y_2}&&\left|
\pk{
\widetilde{M}^{(1)}_n\leq u_n(x_1), \widetilde{M}^{(2)}_n\leq u_n(x_2),
M^{(1)}_n\leq u_n(y_1), M^{(2)}_n\leq u_n(y_2)
} \right.\nonumber\\
&&\qquad \left.-H_{\lambda}^q(x_1,x_2)H_{\lambda}^{1-q}(y_1,y_2)\right|=0,
\EQN
where $H_{\lambda}(x,y)$ is defined in \eqref{add 2} and $a_n,b_n$ are defined in \eqref{eq:ab}.
\item[(2).] Suppose that \eqref{eq:SB} and \eqref{add 3} hold,
then \eqref{eq22}
also holds with $a_n,b_n$ defined in \eqref{anbn B}.
\end{itemize}
\ET
}

\section{Proofs}
\prooftheo{th2.1} By the independence of \Ke{$\vk{S}$ and $(\vk{X}, \vk{Y})$} and \cW{\tE{the} generalised Berman's inequality} (see Theorem 2.1 in \cite{LiShao02} and Lemma 11.1.2 in \cite{leadbetter1983extremes}),
if \eqref{Potter} holds, then
\BQNY
\DUV&=&\pk{-v_i <S_iX_i\leq u_i, 1\le i\le n}-\pk{-v_i<S_iY_i \leq u_i, 1\le i\le n}\\
&\cW{=}&\int_{[0,1]^n}\left(\pk{-\frac{v_i}{s_i}<X_i\leq \frac{u_i}{s_i}, 1\le i\le n}-\pk{-\frac{v_i}{s_i}<Y_i \leq \frac{u_i}{s_i}, 1\le i\le n}\right)
\, d G(s_1)\cdots dG(s_n)\\
&\leq&\frac{2}{\pi}\int_{[0,1]^n}\sum_{1\leq i < j\leq n}\ARIJ
\EXP{-\frac{(w/s_i)^2+(w/s_j)^2}{2(1+\rho_{ij})}} \, dG(s_1)\cdots dG(s_n)\\
&=&\frac{2}{\pi}\sum_{1\leq i < j\leq n} 
\ARIJ
\int_0^1\int_0^1\EXP{-\frac{(w/s)^2+(w/t)^2}{2(1+\rho_{ij})}} \, dG(s) dG(t),
\EQNY
\cE{where \zw{$\rho_{ij}$ and} $\ARIJ$ are defined in \zw{\eqref{ARS}}} \cW{and $w=\min_{1 \le i \le n} \min(\abs{u_i},\abs{v_i})$.}  \cE{Note that}  for $1\leq i, j \leq n$, \wz{$\varepsilon >0$}
\BQNY
&&\int_0^1\EXP{-\frac{1}{2(1+\rho_{ij})}\left(\frac{w}{s}\right)^2} \, dG(s) \\
&\sim&\int_{\frac{1}{\varepsilon+1}}^1\EXP{-\frac{1}{2(1+\rho_{ij})}\left(\frac{w}{s}\right)^2} \, dG(s) \\
&=&\int^{w(1+\varepsilon)}_{w}\pk{S>\frac{w}{s}}\, d\left(1-\EXP{-\frac{1}{2(1+\rho_{ij})}s^2}\right)\\
&=&\int_{0}^{\frac{\varepsilon}{1+\rho_{ij}}w^2}
\pk{S>\frac{w}{w+(1+\rho_{ij})w^{-1}t}}\left(1+\frac{1+\rho_{ij}}{w^2}t\right)
\EXP{-\frac{1}{2(1+\rho_{ij})}(w^2+2(1+\rho_{ij})t+(1+\rho_{ij})^2w^{-2}t^2)}\,dt\\
&\sim&\int_{0}^{\frac{\varepsilon}{1+\rho_{ij}}w^2}
\pk{S>1-\frac{1+\rho_{ij}}{w^2}t}\EXP{-t-\frac{w^2}{2(1+\rho_{ij})}} \, dt\\
&\leq&c_A(1+\rho_{ij})^{\tau}w^{-2\tau}
\EXP{-\frac{w^2}{2(1+\rho_{ij})}}
\int_{0}^{\frac{\varepsilon}{1+\rho_{ij}}w^2}
t^{\tau}\EXP{-t} \, dt\\
&\sim&c_A\Gamma(\tau+1)
(1+\rho_{ij})^{\tau}w^{-2\tau}
\EXP{-\frac{w^2}{2(1+\rho_{ij})}}, \quad \Ke{w\to \IF.}
\EQNY
Consequently, for any $\epsilon>0$ \Ke{and all large $u_i,v_i,i\le n$}
\BQNY
\DUV\leq\frac{2}{\pi}(\Gamma(\tau+1))^2\cW{(c_A^2+\epsilon)}w^{-4\tau}\sum_{1\leq i < j\leq n}
\ARIJ
(1+\rho_{ij})^{2\tau}
\EXP{-\frac{w^2}{1+\rho_{ij}}}.
\EQNY
\wzc{
With similar arguments as above we have
\BQNY
\DUVs&\cW{=}&\int_0^1\left(\pk{-\frac{v_i}{s}<X_i\leq \frac{u_i}{s}, 1\le i \le n}-\pk{-\frac{v_i}{s}<Y_i \leq \frac{u_i}{s}, 1\le i\le n}\right)
\, d G(s)\\
&\le&\frac{2}{\pi}\sum_{1\leq i < j\leq n} 
\ARIJ
\int_0^1\EXP{-\frac{(w/s)^2}{1+\rho_{ij}}} \, dG(s)\\
&\leq&\frac{2^{1-\tau}}{\pi}\Gamma(\tau+1)\cW{(c_A+\epsilon)}w^{-2\rH{\tau}}\sum_{1\leq i < j\leq n}
\ARIJ
(1+\rho_{ij})^{\rH{\tau}}
\EXP{-\frac{w^2}{1+\rho_{ij}}},
\EQNY
}
hence the claim follows.
\QED

\prooftheo{th2.2} According to the independence of the scaling factors with the Gaussian random variables
and the \rH{generalised Berman's inequality} (see Theorem 2.1 in \cite{LiShao02} and Lemma 11.1.2 in \cite{leadbetter1983extremes}) again if \eqref{eq:SB} holds, then  we have
\BQNY
\DUV
&\cW{=}&\int_{[0,\infty]^n}\left(\pk{-\frac{v_i}{s_i}<X_i\leq \frac{u_i}{s_i}, 1\le i \le n}-\pk{-\frac{v_i}{s_i}<Y_i \leq \frac{u_i}{s_i}, 1\le i \le n}\right)
\, d G(s_1)\cdots dG(s_n)\\
&\leq&\frac{2}{\pi}\int_{[0,\infty]^n}\sum_{1\leq i < j\leq n}
\ARIJ
\EXP{-\frac{(w/s_i)^2+(w/s_j)^2}{2(1+\rho_{ij})}} \, dG(s_1)\cdots dG(s_n)\\
&=&\frac{2}{\pi}\sum_{1\leq i < j\leq n}
\ARIJ
\int_0^\infty\int_0^\infty\EXP{-\frac{(w/s)^2+(w/t)^2}{2(1+\rho_{ij})}} \, dG(s) dG(t),
\EQNY
where \zw{$\rho_{ij}$ and} $\ARIJ$ are defined in \zw{\eqref{ARS}}.
Note that for $1\leq i, j \leq n$ and some positive constants $c_1, c_2$, \zc{using similar arguments as in the proof of Theorem 2.1 in \cite{HashorvaWengA}, we have}
\BQNY
&&\int_0^\infty\EXP{-\frac{1}{2(1+\rho_{ij})}\left(\frac{w}{s}\right)^2} \, dG(s)\\
&\sim&\int_{c_1w^{\frac{2}{p+2}}}^{c_2w^{\frac{2}{p+2}}}\EXP{-\frac{1}{2(1+\rho_{ij})}\left(\frac{w}{s}\right)^2} \, dG(s)\\
&\sim&c_BLp\int_{c_1w^{\frac{2}{p+2}}}^{c_2w^{\frac{2}{p+2}}}
s^{\alpha+p-1}\EXP{-Ls^p-\frac{1}{2(1+\rho_{ij})}\left(\frac{w}{s}\right)^2} \, ds\\
&=&c_BLp\left(\frac{w^2}{Lp(1+\rho_{ij})}\right)^{\frac{\alpha+p}{p+2}}
\int_{c_1(Lp(1+\rho_{ij}))^{\frac{1}{p+2}}}^{c_2(Lp(1+\rho_{ij}))^{\frac{1}{p+2}}}
t^{\alpha+p-1}\EXP{-Lp\left(\frac{w^2}{Lp(1+\rho_{ij})}\right)^{\frac{p}{p+2}}\left(p^{-1}t^p+2^{-1}t^{-2}\right)} \, dt\\
&\sim&\sqrt{\frac{2\pi}{p+2}}c_B(Lp)^{\frac{1-\alpha}{p+2}}(1+\rho_{ij})^{\frac{-2\alpha-p}{2(p+2)}}w^{\frac{2\alpha+p}{p+2}}
\EXP{-(1+\rho_{ij})^{-\frac{p}{p+2}}(L^{\frac{2}{p+2}}p^{-\frac{p}{p+2}}+(Lp)^{\frac{2}{p+2}}2^{-1})w^{\frac{2p}{p+2}}}, \quad \zw{w \to \IF.}
\EQNY
Hence for $\epsilon>0$ we have
\BQNY
\DUV\leq
\frac{\cW{4(c_B^2+\epsilon)}(Lp)^{\frac{2(1-\alpha)}{p+2}}}{p+2}w^{\frac{4\alpha+2p}{2+p}}
\sum_{1\leq i < j\leq n}
\ARIJ
(1+\rho_{ij})^{\frac{-2\alpha-p}{p+2}}
\EXP{-2(1+\rho_{ij})^{-\frac{p}{2+p}}Tw^{\frac{2p}{2+p}}},
\EQNY
where $T=L^{\frac{2}{p+2}}p^{-\frac{p}{p+2}}+(Lp)^{\frac{2}{p+2}}2^{-1}$. \wzc{
Proceeding as above
\BQNY
\DUVs&\cW{=}&\int_0^{\IF}\left(\pk{-\frac{v_i}{s}<X_i\leq \frac{u_i}{s}, 1\le i \le n}-\pk{-\frac{v_i}{s}<Y_i \leq \frac{u_i}{s}, 1\le i \le n}\right)
\, d G(s)\\
&\le&\frac{2}{\pi}\sum_{1\leq i < j\leq n}
\ARIJ
\int_0^\infty\EXP{-\frac{(w/s)^2}{1+\rho_{ij}}} \, dG(s)\\
&\leq&2^{\frac{3+2p+\alpha}{2+p}}\pi^{-\frac{1}{2}}\cW{(c_B+\epsilon)}(Lp)^{\frac{1-\alpha}{p+2}}(p+2)^{-\frac{1}{2}}w^{\frac{2\alpha+p}{2+p}}\\
&& \times \sum_{1\leq i < j\leq n}
\ARIJ
(1+\rho_{ij})^{\frac{-2\alpha-p}{2(p+2)}}
\EXP{-(2(1+\rho_{ij})^{-1})^{\frac{p}{2+p}}Tw^{\frac{2p}{2+p}}},
\EQNY}
hence the proof is complete.
\QED

\begin{lemma}\label{le4.1}
Under the conditions of Theorem \ref{th3.1}, for any bounded set $K\subset\{\wz{2,3 \ldots}\}$ we have
$$\lim_{n\to \IF}\pk{S_nX_{n,k}\le u_n, k\in K | S_nX_{n,1}>u_n}=\pk{E/2+\sqrt{\delta_{k-1}}W_k\le \delta_{k-1}, k\in\zw{K}},$$
where $E$ is a standard exponential random variable independent of $\{W_k,k\in K\}$ and the $W_k$
have a jointly Gaussian distribution with mean zero and
$$\E(W_iW_j)=\frac{\delta_{i-1}+\delta_{j-1}-\delta_{|i-j|}}{2\sqrt{\delta_{i-1}\delta_{j-1}}},\quad i,j \in K.$$
\end{lemma}

\prooflem{le4.1} \Ke{A centered} Gaussian random vector
$\bs{X}_n=(X_{n,k},k\in \zw{K\cup\{1\}})^\top$ \zw{with covariance matrix $ B^\top_n B_n=(\wz{\varrho}_{n,|i-j|})_{i,j\in K\cup\{1\}}$}
has stochastic representation
$$(X_{n,k},k\in \zw{K\cup\{1\}})^\top\stackrel{d}{=}R B^\top_n  \Ke{\vk{U}_{m+1}},$$
where $m$ is the cardinality of set $K$, $R$ is a  positive random variable such
that $R^2$ is chi-squared distributed with $ m+1 $ degrees of freedom and independent of $ \Ke{\vk{U}_{m+1}}$
which is a random vector uniformly distributed on the unit sphere of $\mathbb{R}^{ \zw{m+1} }$.
Since $S_n$ \Ke{is} independent of $X_{n,k}$ using Corollary 5 in \cite{MR629795}
we have (set $t_n(y):=u_n+y/u_n$)
$$(S_nX_{n,k},k\in K|S_n X_{n,1}=t_n(y))^{\top}\stackrel{d}{=} R_{\zw{m},y} \hat{B}^\top_n  \Ke{\vk{U}_{\zw{m}}}+t_n(y)\Sigma_{12},$$
where $\Sigma_{12}=(\wz{\varrho}_{n,\zw{k-1}},k\in K)^\top$, \zw{$\hat{B}^\top_n\hat{B}_n=(\varrho_{n,|i-j|}-\varrho_{n,i-1}\varrho_{n,j-1})_{i,j\in K}$}
and $ R_{\zw{m},y} $ is a positive random variable independent of $ \Ke{\vk{U}_m}$ with distribution function $F_{m,y}$ defined by
$$F_{m,y}(x)=
\frac{\int_{t_n(y)}^{((t_n(y))^2+x^2)^{1/2}}\left(s^2-\left(t_n(y)\right)^2\right)^{\frac{ m }{2}-1}s^{1- m } d F_1(s)}
{\int_{t_n(y)}^{\IF}\left(s^2-\left(t_n(y)\right)^2\right)^{\frac{ m }{2}-1}s^{1- m } d F_1(s)}, \quad x>0,
$$
with $F_1$ the distribution function of $\zw{S_n}R$. \Ke{According to Theorem 3.1 in \cite{MR2759158}}
$F_1$ in the Gumbel max-domain of attraction and
\begin{equation}\label{add4.1}
\lim_{n\to \IF}\frac{\pk{S_nX_{n,1}> t_n(y)}}{\pk{S_nX_{n,1}>u_n}}=e^{-y}, \quad \ \ \forall y\inr.
\end{equation}
Hence, by Theorem 3.1 in \cite{MR2247592}
\begin{eqnarray}\label{add4.2}
\Ke{p_{n,y}}&:=&\pk{S_nX_{n,k}\le u_n, k\in K|S_nX_{n,1}=t_n(y)}\nonumber\\
&=&\pk{\frac{u_n(1-\wz{\varrho}_{n,k-1}^2)^{1/2}}{2}Z_{n,k}+\frac{\wz{\varrho}_{n,k-1}}{2}y\le\frac{u_n^2(1-\wz{\varrho}_{n,k-1})}{2}, k\in K}\nonumber\\
&\to&\pk{\sqrt{\delta_{k-1}}W_k+\frac{y}{2}\le \delta_{k-1}, k\in K}, \ \   n\to \IF
\end{eqnarray}
\Ee{uniformly on compact sets of $y$}, where
$$
(Z_{n,k},k\in K)^\top\stackrel{d}{=} R_{m,y} \tilde{B}^\top_n \Ke{\vk{U}_m}, \quad \text{
with }
\tilde{B}^\top_n\tilde{B}_n=\left(\frac{\wz{\varrho}_{n,|i-j|}-\wz{\varrho}_{n,i-1}\wz{\varrho}_{n,j-1}}
{\sqrt{(1-\wz{\varrho}_{n,i-1}^2)(1-\wz{\varrho}_{n,j-1}^2)}}\right)_{i,j\in K}$$
and $\{W_k,k\in K\}$ being jointly Gaussian with zero means and \qE{covariances}
$$\EE{W_iW_j}=\frac{\delta_{i-1}+\delta_{j-1}-\delta_{|i-j|}}{2\sqrt{\delta_{i-1}\delta_{j-1}}}, \ \ 
 i,j\in K.$$
\Ee{Since further}
\begin{eqnarray*}
\pk{S_nX_{n,k}\le u_n, k\in K|S_nX_{n,1}>u_n}=\int_0^{\IF}\Ke{p_{n,y}}d\frac{\pk{S_nX_{n,1}\le t_n(y)}}{\pk{S_nX_{n,1}>u_n}}
\end{eqnarray*}
the proof is established by applying  Lemma 4.4 in \cite{MR2247592} (recall \eqref{add4.1} and \eqref{add4.2}). \QED

\prooftheo{th3.1} According to \zw{\eqref{eq2.4}}, if $1\le k_1 < \ldots <k_s \le n$ \zw{and $k=\min_{1\le i< s}(k_{i+1}-k_{i})$}
then the joint distribution function $F_{k_1,\ldots,k_s}$ of $S_nX_{n,k_1},\ldots,S_nX_{n,k_s}$
satisfies
\begin{eqnarray*}
\left|F_{k_1,\ldots,k_s}(u_n)-\prod_{i=1}^s\pk{S_n X_{\wz{n,k_i}}\le u_n}\right|
\le \mathcal{Q}u_n^{-\wz{2}\tau}n\sum_{\zw{i=k}}^n\frac{|\wz{\varrho}_{n,i}|(1+\zw{|\varrho_{n,i}|})\wz{^{\tau}}}{\sqrt{1-\wz{\varrho}_{n,i}^2}}
\exp\left(-\frac{u_n^2}{1+|\wz{\varrho}_{n,i}|}\right).
\end{eqnarray*}
Suppose now that $1\le i_1 < \ldots< i_p<j_1<\ldots<j_{p^\prime}\le n$ \zw{and $j_1-i_p\ge l_n$}. Identifying $\{k_1,\ldots,k_s\}$
in turn with $\{i_1,\ldots,i_p,j_1,\ldots,j_{p^\prime}\}$, $\{i_1,\ldots,i_p\}$ and $\{j_1,\ldots,j_{p^\prime}\}$,
we thus have
\begin{eqnarray*}\label{add4.3}
|F_{i_1,\ldots,i_p,j_1,\ldots,j_{p^\prime}}(u_n)-F_{i_1,\ldots,i_p}(u_n)F_{j_1,\ldots,j_{p^\prime}}(u_n)|
\le 3\mathcal{Q}u_n^{-\wz{2}\tau}n\sum_{i=\zw{l_n}}^n\frac{|\wz{\varrho}_{n,i}|(1+|\wz{\varrho}_{n,i}|)^{\wz{\tau}}}{\sqrt{1-\wz{\varrho}_{n,i}^2}}
\exp\left(-\frac{u_n^2}{1+|\wz{\varrho}_{n,i}|}\right).
\end{eqnarray*}
By Example 1 in \cite{HASH12} and Table 3.4.4 in
\cite{MR1458613} we have
$$\lim_{n\to\IF}n\pk{S_nX_{n,1}\ge u_n(x)}=e^{-x},\quad \forall\ \ x\inr,$$
where $u_n(x)=a_n x+b_n$ with \wz{$a_n$ and $b_n$ defined in \eqref{anbntau}.}
Consequently, as $n\to \IF$
\begin{equation*}\label{add4.4}
u_n^2(x)=2\ln n-(2\tau+1)\ln\ln n+O(1).
\end{equation*}
Hence, in view of \eqref{add3.1} and \eqref{add3.2},  Theorem 2.1 in \cite{MR877604} implies
$$\limit{n}\left[\pk{\max_{1\le i\le n}S_nX_{n,i}\le u_n(x)}
-\exp\left(-n\pk{S_nX_{n,1}>u_n(x)}\pk{\bigcap_{i=\wz{2}}^{r_n}\{S_nX_{n,i}\le u_n(x)\}|S_nX_{n,1}>u_n(x)}\right)\right] =0.$$
Note that for $m\le j \le r_n$ we have
$$\pk{W>u_n\sqrt{\frac{1-\wz{\varrho}_{n,j}}{1+\wz{\varrho}_{n,j}}}-\frac{y}{u_n}\frac{\wz{\varrho}_{n,j}}{\sqrt{1-\wz{\varrho}_{n,j}^2}}}
\le \mathcal{Q}n^{-\frac{1-\wz{\varrho}_{n,j}}{1+\wz{\varrho}_{n,j}}}\frac{(\ln n)^{\frac{\tau(1-\wz{\varrho}_{n,j})-\wz{\varrho}_{n,j}}{1+\wz{\varrho}_{n,j}}}}{\sqrt{1-\wz{\varrho}_{n,j}^2}},$$
where $W$ is a $N(0,1)$ random variable. The claim can then be established by using similar arguments as in the proof of Theorem 2.1 in \cite{MR1398063} making further use of
\eqref{add3.3} and Lemma \ref{le4.1}. \QED

\COM{
\begin{lemma}\label{le4.2}
Let $\lambda(k,n)=\max(|\lambda_{11}(k,n)|,|\lambda_{12}(k,n)|,|\lambda_{21}(k,n)|,|\lambda_{22}(k,n)|)$.
\COM{
(1)\, If \eqref{eq2.1} holds , we have
$$\lim_{n \to \IF}nw_{\gamma}^{-4\tau}\sum_{k=1}^{n-1}\lambda(k,n)
\EXP{-\frac{w_{\gamma}^2}{1+\lambda(k,n)}} = 0,$$
where
$w_{\gamma}:=w_{\gamma}(z)=a_{n,\gamma}z+b_{n,\gamma},z\inr$ with $a_{n,\gamma}$ and $b_{n,\gamma}$ as
\BQNY\label{anbn gamma}
a_{n,\gamma}=(2\ln n)^{-\frac{1}{2}}, \qquad b_{n,\gamma}=\sqrt{2 \ln n}+(2\ln
n)^{-\frac{1}{2}}\left(\ln \varpi_A-\frac{2\gamma+1}{2}(\ln \ln
n+\ln 2)\right)
\EQNY with $\varpi_A=c_A(2
\pi)^{-\frac{1}{2}}\Gamma(1+\gamma)$.
}
 If \eqref{add 3} holds, then we have
$$\lim_{n\to \IF}nw^{\frac{4\alpha+2p}{2+p}}\sum_{k=1}^{n-1}\lambda(k,n)
\EXP{-2(1+\lambda(k,n))^{-\frac{p}{2+p}}Tw^{\frac{2p}{2+p}}}=
0,$$ where
$w=a_nz+b_n$ with $a_n$ and $b_n$ are defined in \eqref{anbn B}.
\end{lemma}

\prooflem{le4.2}
See Lemma 3.3 in \cite{HashorvaPengWeng}.
\QED

\textbf{Proof.}
Using similar arguments in Lemma 4.3.2 in \cite{leadbetter1983extremes} let $\sigma=\sup_{k\geq1}\rho(k,n)<1$, $\sigma(m,n)=\sup_{k\geq m}\rho(k,n)<1$, $m\leq n$.

(1)\,
Let $\tau$ be a constant
such that $0<\tau <(1-\sigma)/(1+\sigma)$ and $w_{n,\gamma_1,A}=\min\{|u_{n,\gamma_1}(x)|,|u_{n,\gamma_1}(y)|\}$.
According to Lemma \ref{add le3.2} we have
\BQNY
\EXP{-\frac{w_{n,\gamma_1,A}^2}{2}}\sim \mathcal{Q}\frac{(w_{n,\gamma_1,A})^{1+2\gamma_1}}{n} \qquad \qquad
w_{n,\gamma_1,A}\sim \sqrt{2\ln n}.
\EQNY
Split the sum into two parts, i.e.,
\BQNY
n\sum_{k=1}^{n-1}\rho(k,n)
(w_{n,\gamma_1,A})^{-4(\gamma_2-\epsilon)}
\EXP{-\frac{w^2_{n,\gamma_1,A}}{1+\rho(k,n)}}
=n\left(\sum_{k=1}^{[n^\tau]}+\sum_{k=[n^\tau]+1}^{n-1}\right)\rho(k,n)
(w_{n,\gamma_1,A})^{-4(\gamma_2-\epsilon)}\EXP{-\frac{w^2_{n,\gamma_1,A}}{1+\rho(k,n)}}
\EQNY
The first part is dominated by
\BQNY
nn^{\tau}(w_{n,\gamma_1,A})^{-4(\gamma_2-\epsilon)}
\EXP{-\frac{w_{n,\gamma_1,A}^2}{1+\sigma}}
&=&n^{1+\tau}(w_{n,\gamma_1,A})^{-4(\gamma_2-\epsilon)}
\left(\EXP{-\frac{w_{n,\gamma_1,A}^2}{2}}\right)^{\frac{2}{1+\sigma}}\\
&\leq&\mathcal{Q}n^{1+\tau}(w_{n,\gamma_1,A})^{-4(\gamma_2-\epsilon)}
\left(\frac{(w_{n,\gamma_1,A})^{1+2\gamma_1}}{n}\right)^{\frac{2}{1+\sigma}}\\
&\leq&\mathcal{Q}n^{1+\tau-\frac{2}{1+\sigma}}(\ln n)^{\frac{1+2\gamma_1}{1+\sigma}-2(\gamma_2-\epsilon)}.
\EQNY
This tends to zero since $1+\tau-\frac{2}{1+\sigma}<0$. For the second part,
writing $r_n=[n^\tau]$, we have
\BQNY
&&n\sum^{n-1}_{k=r_n+1}\rho(k,n)(w_{n,\gamma_1,A})^{-4(\gamma_2-\epsilon)}\EXP{-\frac{w_{n,\gamma_1,A}^2}{1+\rho(k,n)}}\\
&\leq&n\sigma(r_n,n)(w_{n,\gamma_1,A})^{-4(\gamma_2-\epsilon)}\EXP{-w_{n,\gamma_1,A}^2}\sum^{n-1}_{k=r_n+1}
\EXP{\frac{w_{n,\gamma_1,A}^2\rho(k,n)}{1+\rho(k,n)}}\\
&\leq&n^2\sigma(r_n,n)(w_{n,\gamma_1,A})^{-4(\gamma_2-\epsilon)}\EXP{-w_{n,\gamma_1,A}^2}\EXP{\sigma(r_n,n)w_{n,\gamma_1,A}^2}\\
&\leq&\mathcal{Q}\sigma(r_n,n)(w_{n,\gamma_1,A})^{2+4\gamma_1-4(\gamma_2-\epsilon)}\EXP{\sigma(r_n,n)w_{n,\gamma_1,A}^2}.
\EQNY
Since
\BQNY
\sigma(r_n,n)(w_{n,\gamma_1,A})^{2+4\gamma_1-4(\gamma_2-\epsilon)}&\sim& \sigma(r_n,n)(2\ln n)^{1+\Delta_{\epsilon}}\\
&=&\left(\frac{2}{\tau}\right)^{1+\Delta_{\epsilon}}\sigma(r_n,n)(\ln r_n)^{1+\Delta_{\epsilon}}\\
&\leq&\left(\frac{2}{\tau}\right)^{1+\Delta_{\epsilon}}\sup_{k\geq r_n}\rho(k,n)(\ln k)^{1+\Delta_{\epsilon}} \to 0
\EQNY
and
$$\sigma(r_n,n)(w_{n,\gamma_1,A})^2 \sim 2\sigma(r_n,n)\ln n
\leq \sigma(r_n,n)(2\ln n)^{1+\Delta_{\epsilon}} \to 0$$
as $n \to \infty$, the exponential term above tend to one and the remaining product
tends to zero. Thus we get the first result.

(2)\,  If \eqref{eq:SB} holds, let $\beta$ be a constant
such that $0<\beta <2(1+\sigma)^{-\frac{p}{2+p}}-1$.
According to Lemma \ref{add le3.2} we have
\BQNY
\EXP{-T w_{n,B}^{\frac{2p}{2+p}}}\sim \mathcal{Q}\frac{w_{n,B}^{-\frac{2\alpha}{2+p}}}{n} \qquad \qquad
w_{n,B}\sim \left(\frac{\ln n}{T}\right)^{\frac{2+p}{2p}}.
\EQNY
Split the sum into two parts, i.e.,
\BQNY
&&n\sum_{k=1}^{n-1}\rho(k,n)
(w_{n,B})^{\frac{4\alpha+2p}{2+p}}
\EXP{-2(1+\rho(k,n))^{-\frac{p}{2+p}}T(w_{n,B})^{\frac{2p}{2+p}}}\\
&=&n\left(\sum_{k=1}^{[n^{\beta}]}+\sum_{k=[n^{\beta}]+1}^{n-1}\right)
\rho(k,n)(w_{n,B})^{\frac{4\alpha+2p}{2+p}}
\EXP{-2(1+\rho(k,n))^{-\frac{p}{2+p}}T(w_{n,B})^{\frac{2p}{2+p}}}.
\EQNY
The first part is dominated by
\BQNY
&&nn^{\beta}w_{n,B}^{\frac{4\alpha+2p}{2+p}}
\EXP{-2(1+\sigma)^{-\frac{p}{2+p}}Tw_{n,B}^{\frac{2p}{2+p}}}\\
&=&n^{1+\beta}w_{n,B}^{\frac{4\alpha+2p}{2+p}}
\left(\EXP{-Tw_{n,B}^{\frac{2p}{2+p}}}\right)^{2(1+\sigma)^{-\frac{p}{2+p}}}\\
&\leq&\mathcal{Q}n^{1+\beta}w_{n,B}^{\frac{4\alpha+2p}{2+p}}
\left(\frac{w_{n,B}^{-\frac{2\alpha}{2+p}}}{n}\right)^{2(1+\sigma)^{-\frac{p}{2+p}}}\\
&\leq&\mathcal{Q}n^{1+\beta-2(1+\sigma)^{-\frac{p}{2+p}}} (\ln
n)^{1-\frac{2\alpha}{p}(1-(1+\sigma)^{-\frac{p}{2+p}})}\cH{\to 0}
\EQNY \cH{as $n\to\infty$} since
$1+\beta-2(1+\sigma)^{-\frac{p}{2+p}}<0$.
For the second part, let
$\iota_n=[n^\beta]$, we have \BQNY &&n\sum_{k=\iota_n+1}^{n-1}
\rho(k,n)(w_{n,B})^{\frac{4\alpha+2p}{2+p}}
\EXP{-2(1+\rho(k,n))^{-\frac{p}{2+p}}T(w_{n,B})^{\frac{2p}{2+p}}}\\
&\leq&n^2\sigma(\iota_n,n)
w_{n,B}^{\frac{4\alpha+2p}{2+p}}
\EXP{-2(1+\sigma(\iota_n,n))^{-\frac{p}{2+p}}Tw_{n,B}^{\frac{2p}{2+p}}}\\
&\leq&n^2\sigma(\iota_n,n)
w_{n,B}^{\frac{4\alpha+2p}{2+p}}
\EXP{-2Tw_{n,B}^{\frac{2p}{2+p}}}
\EXP{2T\sigma(\iota_n,n)w_{n,B}^{\frac{2p}{2+p}}}\\
&\leq&\mathcal{Q}\sigma(\iota_n,n)
w_{n,B}^{\frac{2p}{2+p}}
\EXP{2T\sigma(\iota_n,n)w_{n,B}^{\frac{2p}{2+p}}}.
\EQNY
Since
\BQNY
\sigma(\iota_n,n)w_{n,B}^{\frac{2p}{2+p}}\sim T^{-1}\sigma(\iota_n,n)\ln n=\frac{1}{T\beta}\sigma(\iota_n,n)\ln n^{\beta} \to 0
\EQNY
as $n \to \infty$, the exponential term above tend to one and the remaining product
tends to zero. Thus the second result is proved. \QED
}

\wz{
Next, for some index sets $I_n \subset N$ we define
\begin{eqnarray*}
\widehat{\mathbf{M}}(I_n,\bs{\varepsilon}_n):=
\left\{
\begin{array}{ll}
\max\{S_{n,k}\hat{\mathbf{X}}_{n,k}, k \in I_n, \varepsilon_{n,k}=1 \}, & if \ \sum_{k \in I_n}\varepsilon_{n,k}\ge 1;\\
\inf\left\{\mathbf{x}\big{|}\pk{S_{n,k}\hat{\mathbf{X}}_{n,k}\leq\mathbf{x}}>\mathbf{0}\right\}, & otherwise,\\
\end{array}
\right.
\end{eqnarray*}
\begin{eqnarray*}
\widehat{\mathbf{m}}(I_n,\bs{\varepsilon}_n):=
\left\{
\begin{array}{ll}
\min\{S_{n,k}\hat{\mathbf{X}}_{n,k}, k\in I_n, \varepsilon_{n,k}=1 \}, & if \ \sum_{k \in I_n}\varepsilon_{n,k}\ge 1;\\
\inf\left\{\mathbf{x}\big{|}\pk{S_{n,k}\hat{\mathbf{X}}_{n,k}\leq\mathbf{x}}>\mathbf{0}\right\}, & otherwise.\\
\end{array}
\right.
\end{eqnarray*}
For simplicity, we write
$\widehat{\mathbf{M}}_n(\bs{\varepsilon}_n)=\widehat{\mathbf{M}}(\{1,2,\ldots,n\}, \bs{\varepsilon}_n),$
$\widehat{\mathbf{M}}(I_n)=\max\{S_{n,k}\widehat{\mathbf{X}}_{n,k}, k\in I_n\},$
$\widehat{\mathbf{M}}_n=\max\{S_{n,k}\widehat{\mathbf{X}}_{n,k},1\le k\le n\}.$
Similarly we also define $\widehat{\mathbf{m}}_n(\bs{\varepsilon}_n),
\widehat{\mathbf{m}}(I_n), \widehat{\mathbf{m}}_n.$
}

\COM{
Let the independent Gaussian triangular array $\{(\hat{X}_{n,i}^{(1)},\hat{X}_{n,i}^{(2)}),1\le i\le n, n\ge 1\}$ have
stochastic representation
$(\hat{X}_{n,1}^{(1)},\hat{X}_{n,1}^{(2)})
\stackrel{d}{=} (R \cos \theta, R\cos (\theta-\psi_n))$, where
$R$ is a positive random \Ee{variable being independent} of the random variable $\theta$ which is uniformly distributed in $(-\pi,\pi)$ and $\psi_n=\arccos(\lambda_0(n))$ with $\lambda_0(n)$ satisfying \eqref{add 1}. Assume that the distribution function of $SR$ in the Gumbel max-domain
of attraction, where $S$ is a positive random variable and independent of $\{(\hat{X}_{n,i}^{(1)},\hat{X}_{n,i}^{(2)}),1\le i\le n, n\ge 1\}$ . }

\begin{lemma}\label{le4.3}
\zw{Let $\{(\hat{X}_{n,i}^{(1)},\hat{X}_{n,i}^{(2)}),1\le i\le n, n\ge 1\}$ \tE{be} a triangular array of \tE{centered} stationary Gaussian random vectors defined as above with
the correlation $\lambda_0(n)$ satisfying \eqref{add 1} with $\lambda \in (0,\IF)$. Further let $\{S_{n,k},1\le k\le n,n\ge 1\}$ be iid random variables being
independent of $\{(\hat{X}_{n,i}^{(1)},\hat{X}_{n,i}^{(2)}),1\le i\le n, n\ge 1\}$ and satisfying \eqref{eq:SB}.}
Then we have
\BQNY
\lim_{n \to \infty}\pk{-u_n(y_1)<\widehat{m}_n^{(1)}\leq \widehat{M}_n^{(1)}\leq u_n(x_1),-u_n(y_2)< \widehat{m}_n^{(2)}\leq \widehat{M}_n^{(2)}\leq u_n(x_2)}
=H_{\lambda}(x_1,x_2)H_{\lambda}(y_1,y_2).
\EQNY
\end{lemma}

\prooflem{le4.3}
Our proof is similar to that of Theorem 2.1 in \cite{hashorva2012wenglimit}. For any integer $n$ we may write
\BQNY
n\left(1-P(n,x_1,x_2,y_1,y_2)\right)
=nP_1(n,x_1,x_2)+nP_2(n,y_1,y_2)-nP_3(n,x_1,y_2)-nP_4(n,y_1,x_2),
\EQNY
where
\BQNY
&&P(n,x_1,x_2,y_1,y_2):=\pk{-u_n(y_1)<S_{n,1}\hat{X}^{(1)}_{n,1}\leq u_n(x_1),-u_n(y_2)<S_{n,1} \hat{X}^{(2)}_{n,1}\leq u_n(x_2)},\\
&&P_1(n,x_1,x_2):=\pk{S_{n,1}\hat{X}^{(1)}_{n,1}>u_n(x_1)}+\pk{S_{n,1}\hat{X}^{(2)}_{n,1}>u_n(x_2)}-
\pk{S_{n,1}\hat{X}^{(1)}_{n,1}>u_n(x_1),S_{n,1}\hat{X}^{(2)}_{n,1}>u_n(x_2)},\\
&&P_2(n,y_1,y_2):=\pk{S_{n,1}\hat{X}^{(1)}_{n,1}\leq-u_n(y_1)}+\pk{S_{n,1}\hat{X}^{(2)}_{n,1}\leq-u_n(y_2)}-
\pk{S_{n,1}\hat{X}^{(1)}_{n,1}\leq-u_n(y_1),S_{n,1}\hat{X}^{(2)}_{n,1}\leq-u_n(y_2)},\\
&&P_3(n,x_1,y_2):=\pk{S_{n,1}\hat{X}^{(1)}_{n,1}>u_n(x_1),S_{n,1}\hat{X}^{(2)}_{n,1}\leq-u_n(y_2)},\\
&&P_4(n,y_1,x_2):=\pk{S_{n,1}\hat{X}^{(1)}_{n,1}\leq-u_n(y_1),S_{n,1}\hat{X}^{(2)}_{n,1}>u_n(x_2)}.
\EQNY
\zw{The random vector $(\hat{X}_{n,1}^{(1)},\hat{X}_{n,1}^{(2)})$ has the following stochastic representation
$$(\hat{X}_{n,1}^{(1)},\hat{X}_{n,1}^{(2)})
\stackrel{d}{=} (R \cos \theta, R\cos (\theta-\psi_n)),$$ where
$R$ is a positive random \Ee{variable being independent} of the random variable $\theta$ which is uniformly distributed in $(-\pi,\pi)$ and $\psi_n=\arccos(\lambda_0(n))$.
If $S_{n,1}$ satisfy \eqref{eq:SB} and is independent of $(\hat{X}_{n,1}^{(1)},\hat{X}_{n,1}^{(2)})$, using Laplace approximation (see e.g.,\cite{HashorvaWengA}) we have that the
distribution function of $S_{n,1}R$ is in the max-domain of attraction of the Gumbel distribution.
}
Hence,
\zw{according to Remark 2.2 in \cite{HashorvaPengWeng} we have
\begin{equation}\label{addeq1}
\lim_{n\to \IF}n\pk{S_{n,1}\hat{X}_{n,1}^{(1)}>u_n(x)}=e^{-x},\quad x\inr,
\end{equation}
where $u_n(x)=a_nx+b_n$ with $a_n$ and $b_n$ defined in \eqref{anbn B}.
}
Moreover, by Theorem 2.1 in \cite{MR2137118}
\BQNY
&&\lim_{n \to \infty}nP_1(n,x_1,x_2)
=\Phi\left(\lambda+\frac{x_1-x_2}{2\lambda}\right)\exp(-x_2)+
\Phi\left(\lambda+\frac{x_2-x_1}{2\lambda}\right)\exp(-x_1)=:D(x_1,x_2)
\EQNY
and since $(-S_{n,1}\hat{X}^{(1)}_{n,1},-S_{n,1}\hat{X}^{(2)}_{n,1})\stackrel{d}{=}(S_{n,1}\hat{X}^{(1)}_{n,1},S_{n,1}\hat{X}^{(2)}_{n,1})$
\BQNY
&&\lim_{n \to \infty}nP_2(n,y_1,y_2)=D(y_1,y_2).
\EQNY
Since $\lim_{n \to \infty}\lambda_0(n)=1$, $\lim_{n \to \infty }\psi_n=0$ implying
\BQNY
&&\lim_{n \to \infty}nP_3(n,x_1,y_2)\\
&=&\lim_{n \to \infty}n\pk{\wz{S_{n,1}}R\cos(\theta)> u_n(x_1), S_{n,1}R \cos(\theta-\psi_n)\leq-u_n(y_1)}\\
&=&\lim_{n \to \infty}n\pk{S_{n,1}R\cos(\theta)> u_n(x_1), S_{n,1}R \cos(\theta-\psi_n)\leq-u_n(y_1), \cos(\theta)>0, \cos(\theta-\psi_n)<0}\\
&=&\lim_{n \to \infty}n\pk{S_{n,1}R\cos(\theta)> u_n(x_1), S_{n,1}R \cos(\theta-\psi_n)\leq-u_n(y_1), \max\left(-\frac{\pi}{2},-\pi+\psi_n\right)<\theta<-\frac{\pi}{2}+\psi_n}\\
&=&0.
\EQNY
Similarly, we have
$\lim_{n \to \infty}nP_4(n,y_1,x_2)=0.$
Hence for all $x_1,x_2,y_1,y_2\inr$
\BQNY
&&\lim_{n \to \infty}\pk{-u_n(y_1)<\widehat{m}_n^{(1)}\leq \widehat{M}_n^{(1)}\leq u_n(x_1),-u_n(y_2)< \widehat{m}_n^{(2)}\leq \widehat{M}_n^{(2)}\leq u_n(x_2)}\\
&=&\lim_{n \to \infty}(P(n,x_1,x_2,y_1,y_2))^n\\
&=&\lim_{n \to \infty}\left(1-(1-P(n,x_1,x_2,y_1,y_2))\right)^n\\
&=&\lim_{n \to \infty}\left(1-\frac{D(x_1,x_2)+D(y_1,y_2)}{n}+o\left(\frac{1}{n}\right)\right)^n\\
&=&\EXP{-D(x_1,x_2)-D(y_1,y_2)}=H_{\lambda}(x_1,x_2)H_{\lambda}(y_1,y_2),
\EQNY
hence the proof is complete. \QED

\COM{
\begin{remark}
\cW{An important fact for the Gaussian setup is that the distribution function of
$R$ is in the max-domain of attraction of the Gumbel distribution. So 
if $S$ satisfy \eqref{eq:SB}, using Laplace approximation we have
the distribution function of $SR$ in the max-domain of attraction of
the Gumbel distribution.}
\end{remark}
}

\begin{lemma}\label{le4.4}
Under the conditions of Lemma \ref{le4.3}, if the indicator random variables
$\bs{\varepsilon}_n=\{\varepsilon_{n,i},1\le i\le n\}$ are independent of both
$\{(\hat{X}_{n,i}^{(1)},\hat{X}_{n,i}^{(2)}), 1\le i \le n\}$ and $\{S_{n,i},1\le i\le n\}$ and
satisfying condition {\bf E}, then
\BQNY
&&\lim_{n \to \infty}\sup_{x_i, y_i \in \R , i=\{1,2,3,4\}\atop x_1\le x_3, x_2\le x_4, y_1\le y_3, y_2\le y_4 }
\left|\mathbb{P}\left(
-u_n(y_1)< \widehat{m}_n^{(1)}( \bs{\varepsilon}_{n}) \le \widehat{M}_n^{(1)}( \bs{\varepsilon}_{n})\leq u_n(x_1), -u_n(y_2)< \widehat{m}_n^{(2)}(\bs{\varepsilon}_{n}) \le\widehat{M}_n^{(2)}( \bs{\varepsilon}_{n})\leq u_n(x_2),\right.\right.\\
&&\qquad\qquad\qquad\qquad\qquad\qquad\qquad \left.-u_n(y_3)< \widehat{m}_n^{(1)}\le\widehat{M}_n^{(1)}\leq u_n(x_3), -u_n(y_4)< \widehat{m}_n^{(2)} \le \widehat{M}_n^{(2)}\leq u_n(x_4)\right)\\
&&\qquad\qquad\qquad\qquad\qquad\qquad-\left.
\EE{H_{\lambda}^{\eta}(x_1,x_2)H_{\lambda}^{\eta}(y_1,y_2)H_{\lambda}^{1-\eta}(x_3,x_4)H_{\lambda}^{1-\eta}(y_3,y_4)}\right|=0.
\EQNY
\end{lemma}

\prooflem{le4.4}
Using similar arguments as for the derivation of \cite{MR2811020}, let $K_s=\{j: (s-1)\nu+1 \le j \le s\nu\}, 1\le s\le l$, $\nu=[\frac{n}{l}]$,
$\mathbf{x}=(x_1,x_2,x_3,x_4)$, $\mathbf{y}=(y_1,y_2,y_3,y_4)$ and
$\bs{\beta}_n=\{\beta_{n,k}, 1 \le k \le n\}$ be a nonrandom triangular array consisting of $0$'s and $1$'s.
For some random variable $\eta$ such that $0\le \eta \le 1$ a.s., write
\begin{eqnarray*}
B_{\mu,l}=
\left\{\omega: \eta(\omega) \in \Bigg\{
\begin{array}{ll}
[0,\frac{1}{2^l}], & \mu=0,\\
(\frac{\mu}{2^l},\frac{\mu+1}{2^l}], & 0< \mu \le 2^l-1\\
\end{array}
\right\},
\end{eqnarray*}
$$B_{\mu, l, \bs{\beta}_n}=\{\omega: \varepsilon_{n,k}(\omega)=\beta_{n,k}, 1\le k\le n\}\cap B_{\mu, l}.$$
Set
\BQNY
&&P(K_s, \bs{\beta}_n , \mathbf{x}, \mathbf{y})\\
&=&\mathbb{P}\left(
-u_n(y_1)< \widehat{m}^{(1)}(K_s, \bs{\beta}_{n}) \le \widehat{M}^{(1)}(K_s, \bs{\beta}_{n})\leq u_n(x_1), -u_n(y_2)< \widehat{m}^{(2)}(K_s, \bs{\beta}_{n}) \le\widehat{M}^{(2)}(K_s, \bs{\beta}_{n})\leq u_n(x_2),\right.\\
&&\qquad \left.-u_n(y_3)< \widehat{m}^{(1)}(K_s) \le\widehat{M}^{(1)}(K_s)\leq u_n(x_3), -u_n(y_4)< \widehat{m}^{(2)}(K_s) \le \widehat{M}^{(2)}(K_s)\leq u_n(x_4)\right)
\EQNY
and
\BQNY
&&P(n, \bs{\beta}_n , \mathbf{x}, \mathbf{y})\\
&=&\mathbb{P}\left(
-u_n(y_1)< \widehat{m}_n^{(1)}( \bs{\beta}_{n}) \le \widehat{M}_n^{(1)}( \bs{\beta}_{n})\leq u_n(x_1), -u_n(y_2)< \widehat{m}_n^{(2)}(\bs{\beta}_{n}) \le\widehat{M}_n^{(2)}( \bs{\beta}_{n})\leq u_n(x_2),\right.\\
&&\qquad \left.-u_n(y_3)< \widehat{m}_n^{(1)}\le\widehat{M}_n^{(1)}\leq u_n(x_3), -u_n(y_4)< \widehat{m}_n^{(2)} \le \widehat{M}_n^{(2)}\leq u_n(x_4)
\right).
\EQNY
\zw{Using similar arguments as in the proof of Lemma 3.3 in \cite{MR2809873}
 for $n$ large we can choose a positive integer $\tilde{\nu}_n$ such that $l<\tilde{\nu}_n<\nu$ and $\tilde{\nu}_n=o(n)$, by \eqref{addeq1} we have
\begin{eqnarray}\label{addeq2}
&&\left|P(n, \bs{\beta}_n , \mathbf{x}, \mathbf{y})-\prod_{s=1}^lP(K_s, \bs{\beta}_n , \mathbf{x}, \mathbf{y})\right|\nonumber\\
&\le&(4l+2)\tilde{\nu}_n\left(\pk{S_{n,1}\hat{X}_{n,1}^{(1)}\le -u_n(y_1)}+\pk{S_{n,1}\hat{X}_{n,1}^{(1)}> u_n(x_1)}\right.\nonumber\\
&&\qquad \qquad \quad \left.+\pk{S_{n,1}\hat{X}_{n,1}^{(2)}\le -u_n(y_2)}+\pk{S_{n,1}\hat{X}_{n,1}^{(2)}> u_n(x_2)}\right)\nonumber\\
&\to&0, \quad n\to \IF.
\end{eqnarray}
}
Note that
\BQNY
&&1-\frac{\nu\mu}{2^l}\Sigma_1-\nu\left(1-\frac{\mu}{2^l}\right)\Sigma_2+\left(\frac{\sum_{j\in K_s}\beta_{nj}}{\nu}-\frac{\mu}{2^l}\right)\nu(\Sigma_2-\Sigma_1)\\
&\le&P(K_s, \bs{\beta}_n, \mathbf{x}, \mathbf{y})\\
&\le&1-\frac{\nu\mu}{2^l}\Sigma_1-\nu\left(1-\frac{\mu}{2^l}\right)\Sigma_2+\left(\frac{\sum_{j\in K_s}\beta_{nj}}{\nu}-\frac{\mu}{2^l}\right)\nu(\Sigma_2-\Sigma_1)+\nu\Sigma_3,
\EQNY
where
$$\Sigma_1=P_1(n,x_1,x_2)+P_2(n,y_1,y_2)-P_3(n,x_1,y_2)-P_4(n,y_1,x_2),$$
$$\Sigma_2=P_1(n,x_3,x_4)+P_2(n,y_3,y_4)-P_3(n,x_3,y_4)-P_4(n,y_3,x_4)$$
with $P_i(n,z_1,z_2)$'s defined in the proof of Lemma \ref{le4.3}
and
\begin{eqnarray*}
\Sigma_3&=&\sum_{i,j =\{1,2\}}\sum_{t=2}^{\nu}\left(\pk{S_{n,1}\hat{X}_{n,(s-1)\nu+1}^{(i)}>u_n(x_i),S_{n,1}\hat{X}_{n,(s-1)\nu+t}^{(j)}>u_n(x_j)}\right.\\
&&\qquad\qquad+\pk{S_{n,1}\hat{X}_{n,(s-1)\nu+1}^{(i)}>u_n(x_i),S_{n,1}\hat{X}_{n,(s-1)\nu+t}^{(j)}\le -u_n(y_j)}\\
&&\qquad\qquad+\pk{S_{n,1}\hat{X}_{n,(s-1)\nu+1}^{(i)}\le -u_n(y_i),S_{n,1}\hat{X}_{n,(s-1)\nu+t}^{(j)}>u_n(x_j)}\\
&&\qquad\qquad+\left.\pk{S_{n,1}\hat{X}_{n,(s-1)\nu+1}^{(i)}\le -u_n(y_i),S_{n,1}\hat{X}_{n,(s-1)\nu+t}^{(j)}\le -u_n(y_j)}\right).
\end{eqnarray*}
Since $0 \le 1-\frac{\nu\mu}{2^l}\Sigma_1-\nu\left(1-\frac{\mu}{2^l}\right)\Sigma_2 \le 1 $
applying Lemma 3 in \cite{MR2811020} we obtain
\begin{eqnarray}\label{add4.5}
&&\sum_{\mu=0}^{2^l-1}\sum_{\bs{\beta}_n\in\{0,1\}^n}\EE{\left|\prod_{s=1}^lP(K_s,\bs{\beta}_n, \mathbf{x},\mathbf{y})-\prod_{s=1}^l\left[1-\frac{\frac{\mu}{2^l}n\Sigma_1-\left(1-\frac{\mu}{2^l}\right)n\Sigma_2}{l}\right]\right|\II{B_{\mu, l, \bs{\beta}_n}}}\nonumber\\
&\le&\sum_{\mu=0}^{2^l-1}\sum_{\bs{\beta}_n\in\{0,1\}^n}\EE{\sum_{s=1}^l\left|P(K_s,\bs{\beta}_n, \mathbf{x},\mathbf{y})-\left[1-\frac{\frac{\mu}{2^l}n\Sigma_1-\left(1-\frac{\mu}{2^l}\right)n\Sigma_2}{l}\right]\right|\II{B_{\mu, l, \bs{\beta}_n}}}\nonumber\\
&\le&\sum_{\mu=0}^{2^l-1}\sum_{s=1}^l \EE{\frac{\left|\frac{\sum_{j\in K_s}\varepsilon_{n,j}}{\nu}-\frac{\mu}{2^l}\right|}{l}\II{B_{\mu, l}}}n(\Sigma_1-\Sigma_2)+\zw{n} \Sigma_3\nonumber\\
&\le&\sum_{s=1}^l\left[2(2s-1)\left(d\left(\frac{\Xi_{\nu s}}{\nu s}, \eta\right)+d\left(\frac{\Xi_{\nu (s-1)}}{\nu (s-1)}, \eta\right)\right)+\frac{1}{2^l}\right]\frac{n(\Sigma_1-\Sigma_2)}{l}+\zw{n} \Sigma_3,
\end{eqnarray}
where $d(X,Y)$ stands for Ky Fan metric, i.e., $d(X,Y)=\inf\{\varepsilon, \pk{|X-Y|>\varepsilon}<\varepsilon\}$. Furthermore,
\begin{eqnarray}\label{add4.6}
&&\sum_{\mu=0}^{2^l-1}\sum_{\bs{\beta}_n\in\{0,1\}^n}\EE{\left| \prod_{s=1}^l\left[1-\frac{\frac{\mu}{2^l}n\Sigma_1-\left(1-\frac{\mu}{2^l}\right)n\Sigma_2}{l}\right]
-\prod_{s=1}^l\left[1-\frac{\eta n\Sigma_1-\left(1-\eta\right)n\Sigma_2}{l}\right] \right| \I(B_{\mu,l,\bs{\beta}_n})}\nonumber\\
&\le&\sum_{\mu=0}^{2^l-1}\sum_{s=1}^l \EE{ \left|\frac{\mu}{2^l}-\eta\right|\I(B_{\mu,l})}\frac{n}{l}(\Sigma_1+\Sigma_2)\nonumber\\
&\le&\frac{n(\Sigma_1+\Sigma_2)}{2^l}.
\end{eqnarray}
By the fact that  $\lim_{\nu\to \IF}d\left(\frac{\Xi_{\nu s}}{\nu s}, \eta\right)=0$ and utilising \eqref{addeq1}-\eqref{add4.6},
by passing to  limit for $n \to \IF$ and then letting $\nu \to \IF$ we obtain
\BQNY
&&\left|P(n,\vk{\varepsilon}_n,\mathbf{x},\mathbf{y})-\EE{1-\frac{\eta(D(x_1,x_2)+D(y_1,y_2))+(1-\eta)(D(x_3,x_4)+D(y_3,y_4))}{l}}^l\right|\\
&\le&\frac{D(x_1,x_2)+D(y_1,y_2)}{2^{l-1}}\zw{+\frac{1}{l}(e^{-x_1}+e^{-y_1}+e^{-x_2}+e^{-y_2})^2.}
\EQNY
Next, letting $l \to \IF$ implies
\BQNY
\lim_{n \to \infty}\sup_{x_i, y_i \in \R , i=\{1,2,3,4\}\atop x_1\le x_3, x_2\le x_4, y_1\le y_3, y_2\le y_4 }
\left|P(n,\vk{\varepsilon}_n,\mathbf{x},\mathbf{y})
-\EE{H_{\lambda}^{\eta}(x_1,x_2)H_{\lambda}^{\eta}(y_1,y_2)H_{\lambda}^{1-\eta}(x_3,x_4)H_{\lambda}^{1-\eta}(y_3,y_4)}\right|=0,
\EQNY
hence the claim follows.
\QED

\prooftheo{th:m,M}
\COM{Using Corollary \ref{coadd2},
if \wz{\eqref{Potter} and \eqref{eq S gamma}} hold,
we have
\BQNY
SC&:=&\left|\mathbb{P}\left\{
-u_n(y_1)< m^{(1)}_n(\bs{\varepsilon}_n)\le M^{(1)}_n(\bs{\varepsilon}_n)\leq u_n(x_1), -u_n(y_2)< m^{(2)}_n(\bs{\varepsilon}_n)\le M^{(2)}_n(\bs{\varepsilon}_n)\leq u_n(x_2),\right.\right.\\
&&\qquad\left.-u_n(y_3)< m_n^{(1)}\le M^{(1)}_n\leq u_n(x_3),-u_n(y_4)< m_n^{(2)}\le M^{(2)}_n\leq u_n(x_4)
\right\} \\
&& -\mathbb{P}\left\{
-u_n(y_1)< \widehat{m}^{(1)}_n(\bs{\varepsilon}_n) \le \widehat{M}^{(1)}_n(\bs{\varepsilon}_n)\leq u_n(x_1), -u_n(y_2)< \widehat{m}^{(2)}_n(\bs{\varepsilon}_n) \le\widehat{M}^{(2)}_n(\bs{\varepsilon}_n)\leq u_n(x_2),\right.\\
&&\qquad\left. \left.-u_n(y_3)<\widehat{m}^{(1)}_n \le \widehat{M}^{(1)}_n\leq u_n(x_3), -u_n(y_4)< \widehat{m}^{(2)}_n \le\widehat{M}^{(2)}_n\leq u_n(x_4)
\right\}\right|\\
&\leq&\mathcal{Q}nw^{-4\tau}\sum_{i,j=1,2}\sum_{k=1}^{n}|\lambda_{ij}(k,n)|
\EXP{-\frac{w^2}{1+|\lambda_{ij}(k,n)|}}\\
&\leq&\mathcal{Q}nw_{\gamma}^{-4\tau}\sum_{i,j=1,2}\sum_{k=1}^{n}|\lambda_{ij}(k,n)|
\EXP{-\frac{w_{\gamma}^2}{1+|\lambda_{ij}(k,n)|}}.
\EQNY}
If \eqref{eq:SB} holds, \He{by} \zw{\eqref{bb}} \He{for some positive constant $\mathcal{Q}$} we have
\BQNY
&&\left|\mathbb{P}\left(
-u_n(y_1)< m^{(1)}_n(\bs{\varepsilon}_n)\le M^{(1)}_n(\bs{\varepsilon}_n)\leq u_n(x_1), -u_n(y_2)< m^{(2)}_n(\bs{\varepsilon}_n)\le M^{(2)}_n(\bs{\varepsilon}_n)\leq u_n(x_2),\right.\right.\\
&&\qquad\left.-u_n(y_3)< m_n^{(1)}\le M^{(1)}_n\leq u_n(x_3),-u_n(y_4)< m_n^{(2)}\le M^{(2)}_n\leq u_n(x_4)
\right) \\
&& -\mathbb{P}\left(
-u_n(y_1)< \widehat{m}^{(1)}_n(\bs{\varepsilon}_n) \le \widehat{M}^{(1)}_n(\bs{\varepsilon}_n)\leq u_n(x_1), -u_n(y_2)< \widehat{m}^{(2)}_n(\bs{\varepsilon}_n) \le\widehat{M}^{(2)}_n(\bs{\varepsilon}_n)\leq u_n(x_2),\right.\\
&&\qquad\left. \left.-u_n(y_3)<\widehat{m}^{(1)}_n \le \widehat{M}^{(1)}_n\leq u_n(x_3), -u_n(y_4)< \widehat{m}^{(2)}_n \le\widehat{M}^{(2)}_n\leq u_n(x_4)
\right)\right|\\
&\leq&\mathcal{Q}nw^{\frac{4\alpha+2p}{2+p}}\sum_{i,j=1,2}\sum_{k=1}^{n}|\lambda_{ij}(k,n)|
\EXP{-2(1+|\lambda_{ij}(k,n)|)^{-\frac{p}{2+p}}Tw^{\frac{2p}{2+p}}},
\EQNY
\zw{where $w=\min(|u_n(x_i)|,|u_n(y_i)|, 1\le i\le 4)$.}
\He{In view of} Lemma 3.3 in \cite{HashorvaPengWeng}, the sum of the right side of the inequality tends to 0. \He{Thus the claim follows by} Lemma \ref{le4.4}.\QED

{\textbf{Acknowledgments.}}
\qE{We would like to thank the referees for several suggestions which improved our manuscript.}  The authors have been partially support from
the Swiss National Science Foundation grants 200021-140633/1, 200021-134785 and the project RARE -318984 (an FP7 Marie Curie IRSES).

\bibliographystyle{plain}

 \bibliography{BermainINE}
\end{document}